\documentclass[a4paper,12pt]{article}
\usepackage{amsmath}
\usepackage{amsthm}
\usepackage{amssymb}
\usepackage{amscd}
 \makeatletter

 \@addtoreset{equation}{section}
  \makeatother

\theoremstyle{plain}
\newtheorem{thm}{Theorem}[section]
\newtheorem{prop}[thm]{Proposition}
\newtheorem{cor}[thm]{Corollary}
\newtheorem{lem}[thm]{Lemma}
\theoremstyle{definition}

\theoremstyle{remark}
\newtheorem{rem}{Remark}[section]

\newenvironment{namelist}[1]{%
\begin{list}{}
{
\settowidth{\labelwidth}{#1}
\setlength{\leftmargin}{1.1\labelwidth}}
}{%
\end{list}}

\begin{document}
\title{  Prequantization of the moduli space of flat connections over a four-manifold }
\author{ Tosiaki Kori\thanks
{Research supported by
Promotion for Sciences of the Ministry of
 Education and Science in Japan ( no. 16540084 )}\\
Department of Mathematics\\
School of Science and Engineering\\
 Waseda University \\
3-4-1 Okubo, Shinjuku-ku
 Tokyo, Japan.\\ e-mail: kori@waseda.jp
}
\date{ }
\maketitle

\begin{abstract} 
   We introduce a symplectic structure on the space  of connections  in a \(G\)-principal bundle over a four-manifold and the Hamiltonian action on it of the group of gauge transformations which are trivial on the boundary.   The symplectic  reduction becomes the moduli space  of flat connections over the manifold.    On the  moduli space of flat connections we shall construct a hermitian line bundle with connection whose curvature is given by the symplectic form.    This is the Chern-Simons prequantum line bundle.   The group of gauge transformations on the boundary of the base manifold acts on the moduli space of flat connections by an infinitesimally symplectic way.      This action is lifted to the prequantum line bundle by its abelian extension.
\end{abstract}

MSC:   57R; 58E; 81E.

Subj. Class: Global analysis, Quantum field theory.

{\bf Keywords} Chern-Simons prequantization.   Wess-Zumino-Witten actions.     Four-manifolds.   

\medskip

\section{Introduction}

Atiyah and Bott~\cite{AB} showed that the moduli space of flat connections on the trivial \(SU(2)\) bundle over a surface \(\Sigma\) is a compact, finite-dimensional symplectic space.      Ramadas, Singer and Weitsman~\cite{RSW} described the Chern-Simons prequantization of this moduli space.   They showed heuristically how the path-integral quantization of the Chern-Simons action yields holomorphic sections of the prequantum line-bundle when \(\Sigma\) is endowed with a complex structure.     On the other hand, Donaldson proved that if a two-manifold has the boundary then the moduli space of flat connections is a smooth infinite-dimensional symplectic manifold, and it has a  Hamiltonian group action of the gauge transformations on the boundary~\cite{D}.   

In general the moduli space of the \(Lie(G)\)-valued flat connections on a 
a connected compact manifold \(M\) without boundary corresponds bijectively to the conjugate classes of the \(G\)-representations of the fundamental group \(\pi_1(M)\).   
In this paper we study the  Chern-Simons prequantization of the space of flat connections on a four manifold \(M\) generally with non-empty boundary.    

Let \({\cal A}(M)\) be the space of irreducible connections on \(M\)  and let \({\cal G}_0(M)\) be the group of gauge transformations on \(M\) that reduce to the identity transformation when restricted to the boundary \(\partial M\).    We shall prove in section 1 that \({\cal A}(M)\) carries a symplectic structure.    The symplectic form is given by 
\begin{eqnarray*}
\omega_A(a,b)&=&\frac{1}{8\pi^3}\int_MTr[\,(ab-ba)F_A\,]-\frac{1}{24\pi^3}\int_{\partial M}Tr[\,(ab-ba)\,A\,]\\[0.2cm]
&&\qquad +\int_MTr[\,G_A(\ast[a,\ast b])\,]dV\,,
\end{eqnarray*}
for \(a,\,b\in T_A{\cal A}\).     Where \(G_A=(\Delta_A)^{-1}\) is the Green operator of the Dirichlet problem.    The action of \({\cal G}_0(M)\) becomes a Hamiltonian action with the moment map given by the square of curvature \(F_A^2\).   
 Hence the \(0\)-level set of the moment map is the space of flat connections \({\cal A}^{\flat}(M)\).     Let 
\({\cal M}^{\flat}(\partial M)\) be the symplectic reduction, that is, the orbit space of \({\cal A}^{\flat}(M)\) by  \({\cal G}_0(M)\).    We prefer to adopt the notation \({\cal M}^{\flat}(\partial M)\)  than to write \({\cal M}^{\flat}(M)\) because  \({\cal A}^{\flat}(M)/{\cal G}_0(M)\) is near to \({\cal A}^{\flat}(\partial M)\), the space of flat connections over \(\partial M\). 

In section 2 we shall investigate the Chern-Simons prequantization of four manifolds.   By a {\it  prequantization} over \(M\) we mean a  
hermitian line bundle \( {\cal L}(\partial M) \) with connection over \({\cal M}^{\flat}(\partial M)\)  whose curvature is given by the symplectic form of \({\cal M}^{\flat}(\partial M)\),~\cite{ADW,RSW}.     When \(M\) has no boundary \({\cal L}(\partial M)\) becomes a vector space and the {\it Chern-Simons functional} over \(M\) is defined as a vector of this space.   When \(M\) has the boundary the Chern-Simons functional over \(M\) is defined as a section of the pull back line bundle of 
\({\cal L}(\partial M)\) by the quotient map \(p:\,{\cal A}^{\flat}(M)
\longrightarrow {\cal M}^{\flat}(\partial M)\).    The Chern-Simons functional gives a non-vanishing horizontal section over of the pullback line bundle \(p^{\ast}{\cal L}(\partial M)\).

 The set of line bundles \({\cal L}(\partial M)\) and the associated Chern-Simons functionals \({\rm CS}_M\) form a category that represents the operations of union and contraction on the cobordism classes of  base manifolds \(M\), that is, they respect the axioms of Topological Field Theory~\cite{A,K}.     We describe our  {\it Chern-Simons  prequantization} ( CSQ ) as a functor \(CS\) from the category of some classes of four-manifolds with boundary to the category of complex line bundles:
 \[(\,\Gamma=\partial M,\, M\,)  \Longrightarrow (\,{\cal L}(\Gamma),\, CS_M\,) .\]
  As for the axioms of CSQ, we require other than the axioms of TFT two axioms that are characteristic to CSQ.   The first says that 
 \(\pi:\,{\cal L}(\Gamma)
\longrightarrow {\cal M}^{\flat}(\Gamma)\) is a hermitian line bundle with connection whose curvature is equal to the symplectic form on \({\cal M}^{\flat}(\Gamma)\), and the second says that \(\exp 2\pi i\,{\rm CS}_M\) 
is a non-vanishing horizontal section of \(p^{\ast}{\cal L}(\Gamma)\).    

We shall construct the  Chern-Simons  prequantization explicitly for 
\(M=S^4\) and for \(M=D^4\), the 
four-dimensional hemisphere.   Then we give by a functorial method the prequantization over a conformally flat four-manifold with boundary.

The group of pointed gauge transformations \({\cal G}(M)\) acts on 
\({\cal A}^{\flat}(M)\) by infinitesimal symplectic automorphisms.   So is the action of \({\cal G}(\partial M)\) on \({\cal M}^{\flat}(\partial M)\).     We shall discuss in section 3 the lift of this action to the prequantum line bundle \({\cal L}(\partial M)\).    We prove that the abelian extension of   \({\cal G}(\partial M)\) introduced by Mickelsson~\cite{Mi} acts on \({\cal L}(\partial M)\).    
As was discussed in~\cite{K} the geometric description of this abelian extension is given by the four-dimensional Wess-Zumino-Witten model.

\section{ Symplectic structure on the space of connections}

\subsection{Differential calculation on \({\cal A}\)}  

   Let \(M\) be an oriented Riemannian four-manifold with boundary \(\partial M\).   Let \(G=SU(N)\).   The inner product on \(G\) is given by 
\(
<\xi,\eta>=
-Tr(\xi,\eta)\), \( \xi,\,\eta\in Lie(G)=su(N).\)
.   
With this inner product the dual of \(Lie (G)\) is identified with \(Lie (G)\) itself.      
Let \(\pi:\,P\longrightarrow M\) be a principal \(G\)-bundle over \(M\) which is given by a system of transition functions in the Sobolev space \(L^2_s\) for \(s>2\). 
We write \({\cal A}={\cal A}(M)\) for the space of {\it irreducible} \(L^2_{s-1}\) connections, which differ from a smooth connection by an \(L^2_{s-1}\) 1-form  on \(P\) with values in \(Lie (G)\), hence
the tangent space of \({\cal A}\) at \(A\in{\cal A}\) is \(T_A{\cal A}=\Omega^1_{s-1}(M, ad P)\).    Similarly we introduce the space of connections \({\cal A}(\partial M)\) .   It is the \(L^2_{s-\frac12}\) connections that differ from a smooth connection on \(\partial M\) by a \(L^2_{s-\frac12}\) 1-form on \(\partial M\).   
The boundary restriction map  \(r:\,{\cal A}(M)\longrightarrow {\cal A}(\partial M) \) is surjective.   

The curvature of \(A\in {\cal A}\) is 
\[F_A=dA+\frac12[A\wedge A]\in \Omega^2_{s-2}(M, ad P).\]

\(\Omega^3_{s-2}(M, ad P)\), being the dual of \(\Omega^1_{s-1}(M, ad P)=T_A{\cal A}\), is identified with the space of 1-forms on \({\cal A}\).

Here are some differential calculations on \({\cal A}\) that we shall cite from~\cite{BN, DK,S}.

The derivation of a smooth function \(G=G(A)\) on \({\cal A}\) is defined by the functional variation of \(A\):
\begin{equation}
(\delta _A G)a=\lim_{t\longrightarrow 0}\frac{G(A+ta)-G(A)}{t},\quad \mbox{ for \(a\in T_A{\cal A}\)}.
\end{equation}
We have, for example, 
\[(\delta _A A)a=a,\qquad (\delta_AF_A)a=d_Aa.\]
The second follows from the formula
\[F_{A+a}=F_A+d_Aa+a\wedge a.\]

Similarly the derivation of a vector field on \({\cal A}\) or a 1-form on \({\cal A}\)  is  defined as that of a smooth function of \(A\in {\cal A}\) valued in \(\Omega^1(M,ad P)\), respectively,  in \(\Omega^3(M,ad P)\).   
We have, for a vector field \({\bf b}\) and a 1-form \(\beta\), 
\begin{equation}
(\delta_A\,<\beta, {\bf b}>)a=<\beta,\,(\delta_A{\bf b})a>+<(\delta_A\beta)a,{\bf b}>.\end{equation}
where \(<\beta,{\bf b}>\) is the exterior product of \(Lie(G)\) valued forms on \(M\). 

The Lie bracket for vector fields on \({\cal A}\) is seen to have the expression
\begin{equation}
[{\bf a},\,{\bf b}]=(\delta_A{\bf b}){\bf a}-(\delta_A{\bf a}){\bf b}.
\end{equation}

Let \(\widetilde d\) be the exterior derivative on \({\cal A}\).   For a function on \({\cal A}\), \((\widetilde dG)_Aa=(\delta_AG)a\).   
From ( 1.2 ) and ( 1.3 ) we have the following formula for the exterior derivative of a 1-form on \({\cal A}\) :
\begin{eqnarray}
(\widetilde d\theta)_A({\bf a},{\bf b})&=&(\delta_A<\theta,{\bf b}>){\bf a}-(\delta_A<\theta,{\bf a}>){\bf b}-<\theta, [{\bf a},{\bf b}]>\nonumber \\[0.2cm]
&=& <(\delta_A\theta){\bf a},  {\bf b}>-<(\delta_A\theta){\bf b},{\bf a}>.
\end{eqnarray}

Likewise, if \(\omega\) is a 2-form on \({\cal A}\), then 
\begin{equation}
(\widetilde d\omega)_A({\bf a},{\bf b},{\bf c})=(\delta_A\omega({\bf b},{\bf c})){\bf a}+(\delta_A\omega({\bf c},{\bf a})){\bf b}+(\delta_A\omega({\bf c},{\bf a})){\bf b}.\end{equation}

\subsection{ Gauge transformations }

We write \({\cal G}^{\prime}={\cal G}^{\prime}(M)\) for the group of \(L^2_s\) gauge transformations:
\begin{equation}
{\cal G}^{\prime}(M)=\, \Omega^0_s(M, Ad P) .\end{equation}
\({\cal G}^{\prime}\) acts on \({\cal A}\) by 
\begin{equation}
g\cdot A=g^{-1}dg+g^{-1}Ag=A+g^{-1}d_Ag.
\end{equation}
By Sobolev lemma one sees that \({\cal G}^{\prime}\) is a Banach Lie group  and the action is a smooth map of Banach manifolds.   

The group of \(L^2_{s-1/2}\) gauge transformations on the boundary \(\partial M\) is denoted by \({\cal G}^{\prime}(\partial M)\).   We have the restriction map to the boundary:
\[r:\,{\cal G}^{\prime}(M)\longrightarrow {\cal G}^{\prime}(\partial M).\]
Let \({\cal G}_0={\cal G}_0(M)\) be the kernel of the restriction map.   
It is the group of gauge transformations that are identity on the boundary.
\({\cal G}_0\) acts freely on \({\cal A}\) and \({\cal A}/{\cal G}_0\) is therefore a smooth infinite dimensional manifold, while the action of \({\cal G}^{\prime}\) is not free.   
In the following we shall choose \(p_0\in M\) fixed point on the boundary \(\partial M\) and deal with the group of gauge transformations leaving \(p_0\) fixed:
\[{\cal G}={\cal G}(M)
=\{g\in {\cal G}^{\prime}(M);\,g(p_0)=1\,\}.\] 
If \(\partial M=\phi\), \(p_0\) is any point of \(M\) .   \({\cal G}\) act freely on \({\cal A}\).   Let \({\cal A}/{\cal G}\) be the orbit space of this action.   It is a smooth infinite dimensional manifold.    We have  \(Lie({\cal G})=\Omega^0(M, ad P)\).

Correspondingly we have the group \({\cal G}(\partial M)=\{g\in {\cal G}^{\prime}(\partial M);\,g(p_0)=1\,\}\), and the restriction map \(r:\,{\cal G}(M)\longrightarrow {\cal G}(\partial M)\) with the kernel \({\cal G}_0(M)\).   
  We have
\begin{equation}
Lie ({\cal G}_0)=\{\xi\in Lie({\cal G});\,\xi\vert\partial M=0\}.
\end{equation}

The derivative of the action of \({\cal G}\) at \(A\in{\cal A}\) is 
\begin{equation}
d_A=d+[A\wedge\quad]:\, \Omega^0(M, ad P)\longrightarrow \Omega^1(M,ad P).
\end{equation}

Thus the fundamental vector field on \({\cal A}\) corresponding to  
\(\xi \in Lie({\cal G})\) is given by 
\[\xi_{{\cal A}}(A)=\frac{d}{dt}\vert_{t=0}(\exp \,t\xi)\cdot A=d_A\xi.\]
So the tangent space to the orbit at \(A\in{\cal A}\) is 
\begin{equation}
T_A({\cal G}\cdot A)=\{d_A\xi;\,\xi\in \Omega^0(M, ad P)\}
\end{equation}

We have two orbit spaces;
\begin{equation}
 {\cal B}(M)={\cal A}/{\cal G}_0 ,\qquad {\cal C}(M)={\cal A}/{\cal G}.
\end{equation}

\subsection{ Canonical connections on the 
\( {\cal G}_0\)-orbit space and \({\cal G}\)-orbit space.}

Let \({\cal A}={\cal A}(M)\), \({\cal G}={\cal G}(M)\), \({\cal B}={\cal B}(M)\) and \({\cal C}={\cal C}(M)\) be as stated above.   
We shall investigate the horizontal subspaces of the fibrations \({\cal A}\longrightarrow {\cal B}\) and that of  
\({\cal A}\longrightarrow {\cal C}\).   

We have the Stokes formula 
\[\int_{\partial M}f\ast u=\int_M\,d_Af\wedge\ast u- \int_M\,f\ast d^{\ast}_Au,\]
for \(f\in\Omega^0(M,ad P),\, u\in\Omega^1(M,ad P)\).
Hence we see that 
\(d_A\,Lie({\cal G}_0)\) is orthogonal to the space 
\({\rm ker}\,d^{\ast}_A=\{a\in \Omega^1(M,ad P);\, d^{\ast}_Aa=0\}\), 
and .that \(d_A\,Lie({\cal G})\) is orthogonal to \(\{a\in {\rm ker} d^{\ast}_A,\,\, (\ast a)\vert\partial M=0\}\)

Let \(\Delta_A\) be the covariant Laplacian defined as the closed extension of \(d^{\ast}_Ad_A \) with the domain of definition 
\({\cal D}_{\Delta_A}=\{u\in \Omega^0_s(M,ad P);\,u\vert\partial M=0\}\) .   
Since \(A\in{\cal A}\) is irreducible \(\Delta_A:{\cal D}_{\Delta_A}\longrightarrow \Omega_{s-2}^0(M,ad P)\) is an isomorphism.  Let 
\(G_A=(\Delta_A)^{-1}\) be the Green operator of the Dirichlet problem :
\[\left\{\begin{array}{ccc} 
\Delta_A \,u&=& f\\[0.2cm]
u\vert \partial M &=& 0\end{array}\right.\]

\begin{prop}
Let \(A\in{\cal A}\).   
\begin{enumerate}
\item
We have the following orthogonal decomposition:
\begin{equation}
T_A{\cal A}=\{d_A\xi;\,\xi\in Lie({\cal G}_0)\,\}\,\oplus \,H^0_A,\end{equation}
where
\[H^0_A=\{a\in\Omega^1(M, ad P);\, d_A^{\ast}a=0\}.\]
\item
The \({\cal G}_0\)-principal bundle \(\pi:{\cal A}\longrightarrow {\cal B}\) has a natural connection defined by the horizontal subspace \(H^0_A\), 
which is given by the connection form
\(\gamma^0_A=G_Ad_A^{\ast}\).
\item
The curvature form \({\cal F}^0\) of the connection form 
\(\gamma^0\) is given by
\[{\cal F}^0_A(a,b)=G_A(\ast\,[a,\,\ast b])
\qquad\mbox{ for \(a,b \in H^0_A\).}\]
\end{enumerate}
\end{prop}
When \(M\) is compact and \(\partial M=\phi\), the proposition is  well known~\cite{BV, S}., and ours are proved by the same argument.

\vspace{1cm}

\begin{cor}
\begin{equation}
{\cal F}^0_A(a,d_A\xi)=0 \qquad \mbox{for \(\xi\in Lie({\cal G}_0)\)}  .
\end{equation}
\end{cor}

\vspace{1cm}

Now we proceed to the fibration  \({\cal A}\longrightarrow {\cal C}={\cal A}/{\cal G}\).   

Fow a 1-form \(v\), let \(g=K_Av\) denote the solution of the following boundary value problem:
\[\left\{
\begin{array}{ccc} \Delta_A \,g &=& 0\qquad \\[0.2cm]
\ast d_A g\vert\partial M&=&\ast v \vert\partial M.
\end{array}\right.\]

\begin{prop}
Let \(A\in{\cal A}\).   
\begin{enumerate}
\item
We have the orthogonal decomposition:
\begin{equation}
T_A{\cal A}=\{d_A\xi;\,\xi\in Lie({\cal G})\}\,\oplus \,H_A,\end{equation}
where
\[H_A=\{a\in\Omega^1(M, ad P);\, \,d_A^{\ast}a=0,\,\mbox{and } \ast a\vert\partial M=0\}\]
\item
The \({\cal G}\)-principal bundle \(\pi:{\cal A}\longrightarrow {\cal C}\) has a natural connection defined by the horizontal subspace \(H_A\).
\end{enumerate}
\end{prop}

{\it Proof}   

Let \(a\in \Omega^1(M,ad P)\) and \(a=d_A\xi+b\) be the decomposition of 
( 1.11 ), then \(\xi=G_Ad^{\ast}_Aa\) and \(b\in \Omega^1(M,ad P)\), \(d^{\ast}_Ab=0\).   Put \(\eta=K_Ab\).   
Then we have the orthogonal 
decomposition 
\[a=d_A(\xi+\eta)+c,\]
with \(c\in H_A\) and \(\xi+\eta\in Lie({\cal G})\).   
If we write
\begin{equation}
\gamma_A=\gamma^0_A+K_A(I-d_A\gamma^0_A),
\end{equation}
where \(I\) is the identity transformation on \(T_A{\cal A}\), 
then \(\gamma_A\) is a \(Lie({\cal G})-\)valued 1-form which vanishes on \(H_A\) and \(\gamma_A d_A\xi=\xi\),   
that is, \(\gamma_A\) is the connection form.
\hfill\qed

Let \(g=N_Af\) be the solution of Neuman problem:
\[\left\{
\begin{array}{ccc} \Delta_A^{n} \,g &=&f\qquad\\[0.2cm]
\ast d_Ag \vert\partial M&=& 0\qquad\mbox{on \(\partial M\)}.
\end{array}\right.\]
Where \(\Delta_A^{n}\) is the closed extension of \(d^{\ast}_Ad_A \) with the domain of definition 
\({\cal D}_{\Delta^n_A}=\{u\in \Omega^0_s(M,ad P);\,\ast d_Au\vert\partial M=0\}\) .   

\begin{cor}
The curvature form is given by
\begin{equation}
{\cal F}_A(a,b)=N_A(\ast[a,\,\ast b])\qquad\mbox{for \(a,b\in H_A\)}.
\end{equation}
\end{cor}

\vspace{1cm}

\subsection{Symplectic structure on \({\cal A}\)}

For each \(A\in {\cal A}(M)\) we define a sqew-symmetric bilinear form 
on \(T_A{\cal A}\) by:
\begin{eqnarray}
\omega_A(a,b)&=&\omega^0_A(a,b)+\omega^{\prime}_A(a,b)+\phi_A(a,b),\\ [0.2cm]
\omega_A^0(a,b)&=&\frac{1}{8\pi^3}\int_M\,Tr[\,(a\wedge b-b\wedge a)\wedge F_A\,],\\[0.2cm]
\omega^{\prime}_A(a,b)&=&-\frac{1}{24\pi^3}\int_{\partial M}\,Tr[\,(a\wedge b-b\wedge a)\wedge A\,],\\[0.2cm]
\phi_A(a,b)&=&\int_M \,Tr[\,{\cal F}^0_A(a,b)\,]dx,
\end{eqnarray}
for \(a,b\in T_A{\cal A}\).   
Here \({\cal F}^0\) is the curvature form of Proposition 1.1 and \(dx\) is the 
volume form on \(M\).   

Evidently \(\omega\) is non-degenerate.

\begin{thm}
\(\left({\cal A}(M), \omega \right) \) is a symplectic space.
\end{thm}

{\it Proof}

In the following we shall abbreviate \(ab\) for the exterior product \(a\wedge b\).   From the Bianchi's identity we have \(\widetilde d{\cal F}^0=0\), hence \(\widetilde d\phi=0\).   

Differentiating the 2-form \(\omega^0\), we have
\[(\widetilde d\omega^0)_A(a,b,c)=
\delta_A(\omega^0(a,b))(c)+ \delta_A(\omega^0(b,c))(a)+ \delta_A(\omega^0(c,a))(b),\]
for \(a,b,c\in T_A{\cal A}\).   
From the definition we have
\[\delta_A(\omega^0(a,b))(c)=\frac{1}{8\pi^3}\int_M\,Tr[(ab-ba)d_Ac],\]
hence 
\[
(\widetilde d \omega^0)_A(a,b,c)=\frac{1}{8\pi^3}\int_M\,Tr\left [(ab-ba)d_Ac+
(bc-cb)d_Aa+(ca-ac)d_Ab\right].\]
Since
\[d(Tr(ab-ba)c)=
Tr\left [(ab-ba)d_Ac+
(bc-cb)d_Aa+(ca-ac)d_Ab\right],\]
we have 
\[
(\widetilde d \omega^0)_A(a,b,c)=\frac{1}{8\pi^3}\int_Md(Tr(ab-ba) c)=\frac{1}{8\pi^3}
\int_{\partial M}Tr[(ab-ba)c].\]

On the other hand we have
\[(\widetilde d \omega^{\prime})_A(a,b,c)=3\delta_A(\omega^{\prime}(a,b))(c)=-\frac{1}{8\pi^3}\int_{\partial M}Tr[(ab-ba)c].\]
Therefore
\(\widetilde d \omega=0\).

Thus \(({\cal A},\omega)\) is a symplectic manifold.
\hfill \qed

\begin{cor}
\({\cal G}_0\) acts symplectically on \({\cal A}\).
\end{cor}

Note that the action of  \({\cal G}\) is not even infinitesimally symplectic.  

\begin{prop}
The canonical 1-form \(\theta\) that corresponds to the symplectic form \(\omega\) ;
\[ (\widetilde d\theta)_A=\omega_A,
\]
 is given by 
\begin{equation}
\theta_A(a)=-\frac{1}{24\pi^3}\int_M\,Tr[(A F+F A-\frac12 A^3)\, a+ \int_M\,Tr[\gamma_A\,a]dx.
\end{equation}
\end{prop}

\begin{rem}

For \(M\) without boundary the pre-symplectic form \(\omega^0\) and the Hamiltonian action of \({\cal G}={\cal G}_0\) in the next theorem were introduced by Bao and Nair~\cite{BN}.    More generally they gave the pre-symplectic form on n-dimensional manifolds.
\end{rem}

\subsection{ Flat connections }
The space of flat connections are denoted by 
\[{\cal A}^{\flat}(M)=\{A\in{\cal A}(M); F_A=0\},\] 
which we shall often abbreviate to \({\cal A}^{\flat}\).   
The tangent space of \({\cal A}^{\flat}\) is given by
\[T_A{\cal A}^{\flat}=\{a\in \Omega^1(M, ad P); d_Aa=0\}.\]
Here the Sobolev indeces are suppressed.

\begin{thm}~
The action of \({\cal G}_0\) on \({\cal A}\) is a Hamiltonian action and the corresponding moment map is given by 
\begin{eqnarray}
\Phi&:&{\cal A}\longrightarrow (Lie\, {\cal G}_0)^{\ast}=\Omega^4(M, ad P):\,\quad A\longrightarrow F_A^2.\nonumber\\
\nonumber\\
\langle \Phi(A), \xi\rangle&=&\Phi^{\xi}(A)=\frac{1}{8\pi^3}\int_M\,Tr(F_A^2\xi)
\end{eqnarray}
\end{thm}

{\it Proof}

We have 
\[(\delta \Phi^{\xi})_Aa=\frac{1}{8\pi^3}\int_M\,Tr[(d_Aa\wedge F_A+F_A\wedge d_Aa)\xi],\]
and
\begin{eqnarray*}(\delta\Phi^{\xi})_Aa&-&\omega^0_A(a,d_A\xi)=
\frac{1}{8\pi^3}
\int_M\,Tr[(d_Aa\wedge F_A+F_A\wedge d_Aa)\xi-(F_Aa+aF_A)d_A\xi]\\
&&=\frac{1}{8\pi^3}\int_MdTr[(aF_A+F_Aa)\xi]=\frac{1}{8\pi^3}\int_{\partial M}Tr[(aF_A+F_Aa)\xi]=0,
\end{eqnarray*}
 since \(\xi=0\) on \(\partial M\).   On the other hand  we have 
\begin{eqnarray*}
dTr[(Aa-aA)\xi]&=&Tr[(F_Aa+aF_A-Ad_Aa-d_AaA+A^2a+aA^2)\xi]+
Tr[(Aa-aA)d_A\xi]\\
&&=Tr[(ad_A\xi-d_A\xi a)A]\end{eqnarray*}
on \(\partial M\).
Hence
\[
\omega^{\prime}(a,d_A\xi)=-\frac{1}{24\pi^3}\int_{\partial M}
dTr[(Aa-aA)\xi]=0.\]
We have also \(\phi(a,d_A\xi)=0\) from Corollary 1.2.   Therefore
\[ (\delta\Phi^{\xi})_Aa=\omega_A(a,d_A\xi).\]
  
The moment map for the action of \({\cal G}_0\) on \({\cal A}\) is 
\[\Phi:\,A\longrightarrow F_A^2\]
\hfill\qed 

The symplectic quotient of \({\cal A}\) by \({\cal G}_0\) is 
\[{\cal M}^{\flat}=\Phi^{-1}(0)/{\cal G}_0={\cal A}^{\flat}(M)/{\cal G}_0.\]

\begin{rem}~
\({\cal M}^{\flat}\) is the conjugacy classes of \(G\)-representations of \(\pi_1(M,\,\partial M)\).
\end{rem}

\begin{thm}~
\({\cal M}^{\flat}\) is a smooth symplectic manifold.   
The symplectic form on \({\cal M}^{\flat}\) is given by 
\begin{equation}
\omega^{\flat}_{[A]}([a], [b])=\omega^{\prime}_A(a,b)
\end{equation}
for \([A]\in{\cal M}^{\flat}\) and \([a],\,[b] \in T_{[A]}{\cal M}^{\flat}\).
\end{thm}

This is the symplectic reduction theorem of Marsden-Weinstein~\cite{Me} .   
 Let \(({\cal G}_0)_A\) be the isotropy group of \(A\), then the image by 
\((\delta\Phi)_A\) of \(T_A{\cal A}\) is the annihilator of
\[Lie (({\cal G}_0)_A)=\{\xi\in \Omega^0(M, ad P);\,d_A\xi=0,\,
\xi\vert\partial M=0\}=0,\]
since \(A\in{\cal A}\) is irreducible.   
Hence \((\delta\Phi)_A\) is surjective.   In particular \(0\in(Lie\,{\cal G}_0)^{\ast}\) is a regular value.   Thus by the reduction theorem we see that 
\({\cal A}^{\flat}=\Phi^{-1}(0)\) is a co-isotropic submanifold of \({\cal A}\) and the orbit space \({\cal M}^{\flat}={\cal A}^{\flat}/{\cal G}_0\) is a smooth manifold.    The coordinate mappings are described by the implicit function theorem~\cite{DK}.     
For \(A\in{\cal A}^{\flat}\) there is a slice for the \({\cal G}_0\)-action on \({\cal A}^{\flat}\) given by the Coulomb gauge condition:
\[ V_A\subset\{ A+a; \, F_{A+a}=0,\, d^{\ast}_Aa=0\}.\]  
\(V_A\) is a smooth Banach manifold homeomorphic to the following subset of the
tangent space at \(A\):
\[ H_A\cap T_A{\cal A}^{\flat}=\{a\in \Omega^1(M,ad P): \,d_Aa=0,\,d^{\ast}_Aa=0\}.\]
From the definition \(\gamma^0_Aa=0\) for \(a\in H_A\).    
Thus there exists a symplectic structure on \({\cal M}^{\flat}\) given by 
\begin{equation}
\omega^{\flat}_{[A]}([a], [b])=\omega^{\prime}_A(a,b).
\end{equation}
 Here \([A]\in {\cal M}^{\flat}\) denotes the \({\cal G}_0\)-orbit of \(A\in{\cal A}^{\flat}\), and for \([a]\in T_{[A]}{\cal M}^{\flat}\) we take the representative tangent vector to the slice; \(a\in H_A\cap T_A{\cal A}^{\flat}\) .   
\(\omega^{\flat}\) is well defined because we have 
\(g\cdot A=A\) on \(\partial M\) for \(g\in {\cal G}_0\), and \(\omega^{\prime}_A(a,d_A\xi)=0\) for \(\xi\in Lie\,{\cal G}_0\) and \(a\in T_A{\cal A}^{\flat}\).
\hfill\qed 

\begin{cor}
\begin{equation}
T_{[A]}{\cal M}^{\flat}=\{a\in\Omega^1(M.ad P);\,d_Aa=d^{\ast}_Aa=0,\}
\end{equation}
\end{cor}

In the following, when there is a doubt about which manifold is involved,  we shall write 
\({\cal M}^{\flat}(\partial M)\) for the orbit space \({\cal M}^{\flat}={\cal A}^{\flat}(M)/{\cal G}_0(M)\).    

{\bf Example 1}.

For \(M=S^4\), \(\partial S^4=\emptyset\),  we have \({\cal G}_0(S^4)={\cal G}(S^4)\).   
The moduli space of flat connections is 
 \[{\cal M}^{\flat}(\emptyset)={\cal A}^{\flat}(S^4)/{\cal G}(S^4).\]  
Then \({\cal M}^{\flat}(\emptyset)\) is one-point.   

In fact, let \(p_0\in S^4\) and let \(A\in {\cal A}^{\flat}(S^4)\).   Let \(T^A_\gamma(x)\) denote the parallel transformation by \(A\) along the curve \(\gamma\) joining \(p_0\) and \(x\).  
We put 
\(f_A(x)=T^A_\gamma(x)1\in G\).   It is independent of the choice of curve \(\gamma\) joining \(p_0\) and \(x\).   
 Then \(f_A\in {\cal G}(S^4)\).   By the definition  
\(A=df_A\cdot f_A^{-1}\).   

In general \({\cal M}^{\flat}(\emptyset)\) is one-point for a connected and simply connected  compact manifold \(M\) without boundary.

{\bf Example 2.}

For a disc \(D=\{x\in {\rm R}^4;\,\vert x\vert\leq 1\}\) with boundary  \(S^3\), we have \({\cal M}^{\flat}(S^3)\simeq \Omega^3_0G\).   Where 
\[\Omega^3G=\{f\in Map(S^3, G); f(p_0)=1\},\]
and 
\(\Omega^3_0G\) is its connected component of the identity.    
First we note that 
\[{\cal G}(D)\simeq DG=\{f\in Map(D,G);\, f(p_0)=1\},\]
and
\[ {\cal G}_0(D)\simeq D_0G=\{f\in DG;\,f\vert S^3=1 \}.\]
Hence \({\cal G}/{\cal G}_0\simeq \Omega^3_0G\).

As before we put, for \(A\in {\cal A}^{\flat}(D)\),  \(f_A(x)=T^A_{\gamma}(x)1\), \(x\in D\).    We have a well defined bijective map from 
\({\cal A}^{\flat}(D)\) to \(DG\).    In particular, \(f_A=g\) for \(A=dg\,g^{-1}\) with \(g\in {\cal G}_0\) .    
It holds also that \(f_{g\cdot A}(x)=f_A(x)g(x)\) for  \(g\in{\cal G}\) .   
Hence we have the isomorphism
\[ {\cal M}^{\flat}(S^3)= {\cal A}^{\flat}/{\cal G}_0\simeq DG/D_0G\simeq \Omega^3_0G.\]

More generally we have
 \( {\cal M}^{\flat}(\partial M)\simeq {\cal G}/{\cal G}_0\)  
 when \(M\) is a manifold with boundary \(\partial M\) such that \(M\) is made into a simply connected compact manifold by capping a disc to each boundary component.

Now let \(D^{\prime}\subset S^4\) be the complementaly hemisphere of \(D^4 \subset S^4\), the boundary \(\partial D^{\prime}\) is \(S^3\) with the reversed orientation that we denote by \((S^3)^{\prime}=\partial D^{\prime}\).   As above we have 
\[{\cal M}^{\flat}((S^3)^{\prime})={\cal A}^{\flat}(D^{\prime})/{\cal G}_0(D^{\prime}).\]

\begin{lem}~
There is an isomorphism 
\begin{equation}
p:\,{\cal M}^{\flat}(S^3)\longrightarrow {\cal M}^{\flat}((S^3)^{\prime}).
\end{equation}
\end{lem}

{\it Proof}

We define the  map
\begin{equation}
p:\,{\cal A}^{\flat}(D)\longrightarrow {\cal M}^{\flat}((S^3)^{\prime}).
\end{equation}
by
 \[p(A)=[A^{\prime}],\]
where, for \(A\in{\cal A}^{\flat}(D)\) we associated a \(A^{\prime}\in {\cal M}^{\flat}(D^{\prime})\) such that \(A\vert S^3=A^{\prime}\vert S^3\).   The equivalence class \(p(A)=[A^{\prime}]\)  is well defined because, for another \(A^{\prime \prime}\) with 
\(A\vert S^3=A^{\prime \prime}\vert S^3\), there is a \(g^{\prime}\in {\cal G}_0(D^{\prime})\) such that \(A^{\prime \prime}=g^{\prime}\cdot A^{\prime}\).   If \(p(A)=p(B)=[A^{\prime}]\)  then \(A\vert S^3=A^{\prime}\vert S^3=B\vert S^3\), so there is a \(g\in {\cal G}_0(D)\) such that \(B=g\cdot A\).   Hence the kernel of \(p\) is \({\cal G}_0(D)\).    Thus \(p\) induces an  isomorphism:
\begin{equation}
p:\,{\cal M}^{\flat}(S^3)\longrightarrow {\cal M}^{\flat}((S^3 )^{\prime}).
\end{equation}
\hfill\qed 

 The boundary restriction map \(r:\,{\cal A}(M)\longrightarrow {\cal A}(\partial M)\) induces a map
\[\Bar{r}:{\cal M}^{\flat}(\partial M)={\cal A}^{\flat}(M)/{\cal G}_0\longrightarrow {\cal A}^{\flat}(\partial M).\]
If \(A_1,\,A_2\in{\cal A}^{\flat}(M)\) have the same boundary restriction, then the corresponding holonomy groups coincide on the boundary, and defines a 
\(g\in{\cal G}_0\) such that \(g\cdot A_1=A_2\).   Therefore \(\Bar{r}\) is a injective map.

\subsection{The action of \({\cal G}(M)\)}

By the action of the group of gauge transformations  \({\cal G}\)  on \({\cal A}^{\flat}\) we have the orbit space \({\cal N}^{\flat}={\cal A}^{\flat}/{\cal G}\).    Then we have a fibration \({\cal M}^{\flat} \longrightarrow {\cal N}^{\flat}\) with the fiber \({\cal G}(\partial M)={\cal G}/{\cal G}_0\).
We note the fact that any vector which is tangent to the \({\cal G}\)-orbit through \(A\in{\cal A}^{\flat}\) is in \(T_A{\cal A}^{\flat}\).   
 Propositions 1.3 yields the following Proposition.

\begin{prop}~
Let \(A \in{\cal M}^{\flat}\).   
\begin{enumerate}
\item
We have the following decomposition 
\begin{eqnarray}
T_{A}{\cal M}^{\flat}&=&\{d_A\xi;\,\xi\in Lie({\cal G}(\partial M))\}\,\oplus \,H_A^{\flat},\\[0.2cm]
\mbox{where}\nonumber\\
H_A^{\flat}&=&\{a\in\Omega^1(M, ad P);\, \,d_Aa=d_A^{\ast}a=0,\,\mbox{and } \ast a\vert\partial M=0\}.
\nonumber
\end{eqnarray}
\item
The \({\cal G}(\partial M)\)-principal bundle \(\pi:{\cal M}^{\flat}\longrightarrow {\cal N}^{\flat}\) has a natural connection defined by the horizontal space \(H_A^{\flat}\).   The connection form is given by \(K_A(I-d_A\gamma^0_A)\).

\end{enumerate}
\end{prop}

\vspace{0.5cm}

The action of \({\cal G}\) on \({\cal A}\) is far from symplectic.   
But on \({\cal A}^{\flat}\) it acts infinitesimally symplectic.   
In fact, we have \(d_A\xi\in T_A{\cal A}^{\flat}\) for \(\xi \in Lie({\cal G})\) and \(A\in {\cal A}^{\flat}\), and if  we 
denote by \(L_\xi\) the Lie derivative by the fundamental vector field corresponding to \(\xi\in Lie\,{\cal G}\) we have;
\begin{eqnarray*}
(L_\xi\omega)_A(a,b)&=&(\widetilde d\,i_{d_A\xi}\omega)_A(a,b)
=\delta_A\,(i_{d_A\xi}\omega_A(b))(a)-\delta_A\,(i_{d_A\xi}\omega_A(a))(b)\\
&=&-\frac{1}{24\pi^3}\int_{\partial M}\,Tr[bd_A\xi-d_A\xi b)a]+
\frac{1}{24\pi^3}\int_{\partial M}\,Tr[ad_A\xi-d_A\xi a)b]\\
&=&-\frac{1}{12\pi^3}\int_{\partial M}\,Tr[(ab-ba)d_A\xi]
= -\frac{1}{12\pi^3}\int_{\partial M}\,d\,Tr[(ab-ba)\xi]=0,
\end{eqnarray*}
for  \(A\in{\cal A}^{\flat}\) and for 
\(a,b\in T_A{\cal A}^{\flat}\).   

\({\cal G}_0\) being a normal subgroup of \({\cal G}\) the action of 
the quotient group \({\cal G}/{\cal G}_0\)  on
 \( ({\cal M}^{\flat}, \omega^{\flat})\) is also 
infinitesimally symplectic.      

{\bf Example }

The same argument as in Examples 1 and 2 by using paralel transformations along the curves in \(S^3\) yields that 
\[\Omega^3G\simeq {\cal A}^{\flat}(S^3).\]    
So we have an injective mapping  \({\cal M}^{\flat}(S^3)\longrightarrow  {\cal A}^{\flat}(S^3)\).   It corresponds to the embedding \(\Omega^3_0G\longrightarrow \Omega^3G\).     Since \({\cal M}^{\flat}(S^3)\) is a fiber bundle over \({\cal N}^{\flat}(S^3)\) with the fiber \({\cal G}/{\cal G}_0\simeq \Omega^3_0G\), 
\[{\cal N}^{\flat}(S^3)=\,\mbox{ one point}.\]

\vspace{1cm}

\section{ Chern-Simons prequantization over four-manifolds }

Let \(M\) be an oriented Riemannian four-manifold with boundary and \(P=M\times SU(N)\) be the trivial \(SU(N)-\) principal bundle.   
In this section we shall deal with the prequantization of the moduli space of flat connections, that is, over the symplectic manifold \(({\cal M}^{\flat}(\partial M),\omega^{\flat})\) we shall construct  a line bundle with connection whose curvature is 
\(\omega^{\flat}\).     

\subsection{ Descent equations}

Let \(\pi:\,P\longrightarrow M\) be a principal \(G\)-bundle over \(M\).   Let \(\Omega^q\) be the differential \(q\)-forms on \(P\) and let \(V^q\) be the vector space of polynomials \(\Phi(A)\) of \(A\in {\cal A}\)  and its curvature \(F_A\) that take values in \(\Omega^q\).   The group of gauge transformations \({\cal G}\) acts on \(V^q\) by \((g\cdot \Phi)(A)=\Phi(g^{-1}\cdot A)\).   We shall investigate the double complex 
\[C^{p,q}=C^p({\cal G}, V^{q+3}),\] 
that is doubly graded by the chain degree \(p\) and the differential form
degree \(q\).     The coboundary operator \(\delta: C^p\longrightarrow C^{p+1}\) is given by
\begin{eqnarray*}
(\delta\, c^p)(g_1,g_2,\cdots,g_{p+1})&=&
g_1\cdot c^p(g_2,\cdots,g_{p+1})+(-1)^{p+1}c^p(g_1,g_2,\cdots,g_p)\nonumber\\[0.5cm]
&+&\sum_{k=1}^p(-1)^kc^p(g_1,\cdots,g_{k-1},g_kg_{k+1},g_{k+2},\cdots,g_{p+1}).
\end{eqnarray*}
The following Proposition is a precise version of the Zumino's descent equation.     Though stated for \(n=3\) it holds for general \(n\).   The first equation is nothing but the Zumino's equation~\cite{Z}.     The author learned the second equation with its transparent proof from Y. Terashima~\cite{T}.    For \(n=3\) stated here the calculation was appeared in~\cite{Mi,Mic}.

\begin{prop}~
Let \(c^{0,3}\in C^{0,3}\) be defined by \[c^{0,3}(A)=Tr\,F_A^3, \qquad A\in {\cal A}.\]
Then there is a sequence of cochains \(c^{p,q}\in C^{p,q}\), \(0\leq p,q\leq 3\)  that satisfies the following relations:
\begin{eqnarray}
dc^{p,3-p}+(-1)^p\delta c^{p-1,3-p+1}&=&0\\
dc^{p,2-p}+(-1)^{p-1}\delta c^{p-1,3-p}&=&c^{p,3-p}\\
c^{p,q}=0,\qquad \mbox{if}\,p+q\neq 2,3\nonumber
\end{eqnarray}
\end{prop}
   Each term is given by the following form:
\begin{eqnarray*}
c^{0,2}(A)&=&Tr\,(AF^2-\frac12
A^3F+\frac1{10}A^5).\\
c^{1,2}(g)&=& \frac1{10}
Tr(dg\cdot g^{-1})^5\\
c^{1,1}(g; A)&=& = Tr[\,-\frac12V(AF+FA-A^3)+
\frac14(VA)^2+\frac12V^3A]\\
&&\qquad\mbox{where \(V=dg\,g^{-1}\)}\\
c^{2,1}(g_1,g_2)&=&c^{1,1}(g_2 ;\,g_1^{-1}dg_1\,)\\
c^{2,0}(g_1,g_2;A)&=&-Tr\lbrack \frac12(dg_2g_2^{-1})(g_1^{-1}dg_1)
(g_1^{-1}Ag_1)-\frac12(dg_2g_2^{-1})(g_1^{-1}Ag_1)
(g_1^{-1}dg_1)\rbrack .
\end{eqnarray*}

\begin{rem} 
We write here the two-dimensional analogy which is familiar to the readers.   In this case we consider the double complex 
\(C^{p,q}=C^p({\cal G},V^{q+2})\).   
Put 
\(c^{0,2}=Tr F^2\in C^{0,2}\), then there is a sequence of cochains \( c^{p,q}\in C^{p,q}\), \(0\leq p,q \leq 2 \) that satisfies the following relations:
\begin{eqnarray*}
dc^{0,1}=c^{0,2},\quad \delta c^{1,0}&=&c^{2,0}\\
dc^{1,0}+\delta c^{0,1}&=&c^{1,1}\\
dc^{2,0}-\delta c^{1,1}&=&0\\
c^{p,q}=0,\qquad \mbox{if}\,p+q\neq 2,1\nonumber
\end{eqnarray*}
Each terms are given as follows.
\begin{eqnarray*}
c^{0,1}(A)&=&Tr\,(AF-\frac13A^3).\\
c^{1,0}(A; g)&=& Tr(dg\,g^{-1}\,A)\\
c^{1,1}(g)&=&\frac13Tr(dg\,g^{-1})^3\\
c^{2,0}(g_1,g_2)&=&c^{1,0}(g_1^{-1}dg_1;\,dg_2\,g_2^{-1})
\end{eqnarray*}
We note that the relation 
\(dc^{2,0}-\delta c^{1,1}=0\) represents the Polyakov-Wiegmann formula appeared in~\cite{PW},  and the relation 
\(\delta c^{2,0}=0\) is the cocycle condition for the central extension of the loopgroup \(LG\).
\end{rem}

\subsection{ Prequantization over \(\emptyset=\partial S^4\)}

\subsubsection{4-dimensional Polyakov-Wiegmann formula}

Let \(G=SU(N)\).    We assume in the following that \(N\geq 3\).   
Then \(\pi_4(G)=1\) and \(\pi_5(G)\simeq{\rm Z}\).   
Let \(P\) be the trivial \(G\)-bundle over \(S^4\).    
Let \({\cal A}\) be the space of irreducible connections on \(P\).    Let \({\cal G}(M)\) be the group of based gauge transformations.     We have \({\cal G}=S^4G\), the set of smooth mappings from \(S^4\) to \(G\) based at some point.

The Chern-Simons form on \(P\) is defined by 
\begin{equation}
c^{0,2}(A)=Tr\,(AF^2-\frac12A^3F+\frac1{10}A^5).
\end{equation}
 We have then
\(Tr(F^3)=d\,c^{0,2}(A)\).

The variation of the Chern-Simons form along the \({\cal G}\)-orbit is described by the equation ( 2.2 ):
\begin{equation}
c^{0,2}(g\cdot A)-c^{0,2}(A)=d\,c^{1,1}(g; A)+c^{1,2}(g).
\end{equation}
   
Let \(D^5\) be 
the five dimensional ball with boundary \(S^4\).    
Then \(g\in S^4 G\) is extended to 
\(D^5G\), in fact, we have such an extension by virtue of    
\(\pi_4(G)=1\).  

Put
\begin{eqnarray}
\Gamma(g;A)&=&\frac{i}{24\pi^3}\int_{S^4}c^{1,1}(g;A)+C_5(g),
\nonumber\\[0.5cm]
C_5(g)&=&\frac{i}{24\pi^3}\int_{D^5}c^{1,2}(g)\\
&=& \frac{i}{240\pi^3}\int_{D^5}Tr(dg\cdot g^{-1})^5.
\end{eqnarray}
\(C_5(g)\) may depend on the extension,
 but since \(\pi_5(G)={\rm Z}\) the difference 
of two extensions is an integer, and
\(\exp 2\pi i C_5(g)\,\)is independent of the
extension.   

\begin{lem}[Polyakov-Wiegmann]~\cite{K,Mi,PW}~  For \(f,\,g\in
S^4 G\) we have
\begin{equation}
C_5(fg)=C_5(f)+C_5(g)+\gamma(f,g)\quad
\mod{\bf Z},\end{equation}
where
\begin{eqnarray*}
\gamma(f,g)&=&\frac{i}{24\pi^3}\int_{S^4}\,c^{2,1}(f,g)\\
\\
&=&\frac{i}{48\pi^3}\int_{S^4}
Tr\lbrack (dgg^{-1})(f^{-1}df)^3+\frac12(dgg^{-1}f^{-1}df)^2+
\\
\\
&&\mbox{  }\qquad\mbox{ }+(dgg^{-1})^3(f^{-1}df)\rbrack.
\end{eqnarray*}
\end{lem}
{\it proof}

From (2,2)  we have
\[\delta c^{1,2}+dc^{2,1}=0.\]
Integrating it over \(D^5\) yields the desired equation.
\hfill\qed 

\begin{lem}
\begin{equation}
\Gamma(fg,A)=
\Gamma(g,\,f\cdot A)+\Gamma(f,A)
\end{equation}
\end{lem}
{\it proof}

From (2,2)  we have
\[\delta c^{0,2}+dc^{1,1}=c^{1,2}.\]
Hence \(\delta dc^{1,1}=dc^{2,1}\).   
Integrating it over \(D^5\), the Polyakov-Wiegmann formula yields the result.
\hfill\qed 

\subsubsection{Prequantum line bundle over \(\emptyset=\partial S^4\)}

We consider the \(U(1)-\)valued function on \({\cal A}\times {\cal G}\):
\begin{equation}
\Theta(g,A)=\exp 2\pi i\Gamma(g,A).
\end{equation}
Lemma 2.3 yield 
the following cocycle condition:
\begin{equation}
\Theta(g,A)\Theta(h,g\cdot A)=\Theta(gh,A).
\end{equation}
Therefore if we define the action of \({\cal G}(S^4)\) on \({\cal A}(S^4)\times {\rm C}\) 
by 
\[\left(g, \,(A,c)\,\right)\longrightarrow \left(g\cdot A,\,\Theta(g,A)c\right),\]
we have a complex line bundle:
\begin{equation}
{\cal L}={\cal A}\times {\rm C}/{\cal G}\longrightarrow {\cal B}={\cal A}/{\cal G}.
\end{equation}
\(\Theta\) being \(U(1)\)-valued, it is a hermitian line bundle.

We restrict \({\cal L}\) to \({\cal M}^{\flat}(\emptyset)={\cal A}^{\flat}/{\cal G}\) and have the line bundle 
\begin{equation}
{\cal L}(\emptyset)={\cal A}^{\flat}\times {\rm C}/{\cal G}\,\longrightarrow \,{\cal M}^{\flat}(\emptyset).
\end{equation}
We call it the {\it prequantum line bundle} over \(\emptyset\); an empty three-manifold.

Since \({\cal M}^{\flat}(\emptyset)\) is one point \({\cal L}(\emptyset)\)  is isomorphic to the complex line \({\rm C}\). 
\vspace{0.5cm}

\subsubsection{Chern-Simons functional over \(S^4\)}
Now we introduce the 
{\it Chern-Simons functional } \({\rm CS}_{S^4}\).  

We put, for  \(A\in{\cal A}^{\flat}(S^4)\),
\begin{equation}
{\rm CS}_{S^4}(A)=\frac{i}{24\pi^3}\int_{D^5}c^{0,2}(\widetilde A)=\frac{i}{240\pi^3}\int_{D^5}\,Tr(\widetilde A^5),
\end{equation}
where \(\widetilde A\) is the extension of \(A\) to \(D^5\).
In general this integral may depend on the extension  \(\widetilde A\), but, 
as we saw in 1.5, a flat connection \(A\) on \(S^4\) is expressed as \(A=df\,f^{-1}\) by  \(f=f_A\in {\cal G}(S^4)\), so we have
\[c^{0,2}(A)=\frac1{10}Tr(df\cdot f^{-1})^5=c^{1,2}(f).\]
Hence, if we extend \(A=df\,f^{-1}\) to \(\widetilde A=d\widetilde f\,\widetilde f^{\,-1}\) on \(D^5\), we have 
\[{\rm CS}_{S^4}(A)=\frac{i}{24\pi^3}\int_{D^5}c^{0,2}(\widetilde A)= 
\frac{i}{24\pi^3}\int_{D^5}c^{1,2}(\widetilde f)=C_5(f),\]
by virtue of the fact \(\pi_5(G)={\rm Z}\) the last quantity is defined independently of the extension of \(f\) to \(D^5\) modulo an integer, .        Thus \(\exp 2\pi i\,{\rm CS}_{S^4}(A)\) for 
\(A\in{\cal A}^{\flat}(S^4)\) is well defined.     

From ( 2.7 ) we have
\[{\rm CS}_{S^4}(g\cdot A)=\Gamma(g, A)+{\rm CS}_{S^4}(A),\qquad\mbox{ mod. \({\rm Z}\),}\]
which implies
\begin{equation}
\exp 2\pi i\,{\rm CS}_{S^4}(g\cdot A)=\Theta(g,A)\exp 2\pi i\,{\rm CS}_{S^4}(A).
\end{equation}
This formula defines a non-vanishing element 
\([A,\,\exp 2\pi i\,{\rm CS}_{S^4}(A)]\in 
{\cal L}(\emptyset)\).   
Hence the Chern-Simons functional \({\rm CS}_{S^4}\) gives the trivialization.
\begin{eqnarray}
 {\cal L}(\emptyset)&\,\simeq \,&{\rm C},\\[0.2cm]
[A,c\exp 2\pi i\,{\rm CS}_{S^4}(A)]&\longrightarrow& c.\nonumber
\end{eqnarray}

\begin{lem}
\begin{equation}
\left(\widetilde d\, {\rm CS}_{S^4}\right)_A(a)\,=\frac{i}{48\pi^3}\int_{S^4}\,Tr[ A^3\,a],
\end{equation}
for any \(A\in{\cal A}^{\flat}(S^4)\).
\end{lem}

{\it  proof}

The lemma follows from the following calculations.   
For any \(A\in{\cal A}(S^4)\) and for any extension \(\widetilde A\) of \(A\) to \(D^5\), we have
\[\left(\widetilde d \,Tr[\widetilde A^{\,5}]\right)(a)=5Tr[\widetilde A^{\,4}\,a],\]
and 
\[d\,Tr[A^3\,a]=Tr[ \,\widetilde A^{\,4}\,a\,],\]
where \(a\in T_{\tilde A}{\cal A}(D^5)\).   
And we have
\[\left(\widetilde d\, {\rm CS}_{S^4}\right)_A(a)=\frac{i}{48\pi^3}\int_{D^5}\,Tr[\widetilde A^{\,4}a]=\frac{i}{48\pi^3}\int_{S^4}\,Tr[A^3\,a]
\]
\hfill\qed

We see from this lemma that if we define the connection form \(\theta\) on the line bundle \({\cal L}(\emptyset)\) by \(-2\pi i\) times the right hand side of ( 2.16 ) then 
\(\exp 2\pi i\,{\rm CS}_{S^4}\) is a horizontal section of the pull-back line bundle \(p^{\ast}{\cal L}(\emptyset)\) by \(p:\,{\cal A}^{\flat}(S^4)\longrightarrow {\cal M}^{\flat}(\emptyset)\).     

\vspace{1cm}

\subsection{Prequantization over 4-dimensional discs }

\subsubsection{Pre-quantum line bundles over \(S^3\) and \((S^3)'\) }

Let \(\Omega^3G\) be the set of smooth mappings from
\(S^3\) to \(G=SU(N)\) that are based at some point.   \(\Omega^3G\) is not connected but divided into  the
connected components by \(\deg\).   We put 
\begin{equation}\Omega^3_0G=\{g\in \Omega^3G;\, \deg g=0\}.\end{equation}

The oriented 4-dimensional disc with boundary \(S^3\) is denoted by \(D\), while that with oposite orientation  is denoted by \(D^{\prime} \).   The composite cobordism of \(D\) and \(D^{\prime}\) becomes \(S^4\).   We write  \(DG=Map(D,G)\), the set of smooth mappings from \(D \) to \(G\) based at a \(p_0\in \partial D=S^3\), and  similarly \(D^{\prime} G=Map(D^{\prime}, G)\).   The restriction to \(S^3\) of a \(f\in DG\) has degree 0; \(\,f\vert S^3\in \Omega^3_0G\).    We denote \(D_0G=\{f\in DG;\,f\vert S^3=1\}\) and  \(D^{\prime}_0G=\{f^{\prime}\in D^{\prime}G;\,f^{\prime}\vert S^3=1\}\).    
 The upper prime will indicate  that the function expressed by
the letter  is defined on \(D^{\prime} \), for example, \(1^{\prime}\) is the 
constant function \(D^{\prime}\ni x\longrightarrow 1'(x)=1\in G\), while 
 \(1\) is the constant function 
\(\, D\ni x\longrightarrow 1(x)=1\in
G\).    For  \(g\in DG \) and  \(g^{\prime}\in D^{\prime}G \) such that \(g\vert S^3=g^{\prime}\vert S^3\), we write by 
\(g\vee g^{\prime}\) the gauge transformation on \(S^4\) that are obtained by sewing \(g^{\prime}\) and \(g\).   

Let \(P\) be the trivial \(G=SU(N)\)-principal bundle over \(D\) with \(N\geq 3\).   
 For the groups of gauge transformations on \(D\),   
we have \({\cal G}(D)=DG\),  \({\cal G}_0(D)=D_0G\) and \({\cal G}(S^3)=\Omega^3_0G\).   Similarly for the trivial \(G\)-bundle over \(D'\) we have  \({\cal G}(D^{\prime})=D^{\prime}G\) and  \({\cal G}_0(D')=D^{\prime}_0G\).  

Let \(\left(\,\,{\cal M}^{\flat}(S^3)={\cal A}^{\flat}(D)/D_0G \,,\,\omega^{\flat}\,\right)\) be the moduli space of flat connections over \(D\) with the symplectic structure \(\omega^{\flat}\).   
 We shall construct a prequantum line bundle over \({\cal M}^{\flat}(S^3)\).    Let \(\tilde P\) be the  trivial extension of \(P\) to \(S^4\).     
We extend any \(g\in {\cal G}_0(D)\) across the boundary \(S^3\) by defining it to be the identity transformation on \(\tilde P\vert D^{\prime}\).
Then \(g\vee 1^{\prime}\in S^4G\).   We put, for  \(A\in{\cal A}^{\flat}(D)\) and \(g\in D_0G\),  
\begin{eqnarray}
\Gamma_D(g;A)&=&\frac{i}{24\pi^3}\int_{D}c^{1,1}(g;A)+C_5(g\vee 1^{\prime}),.\\
\nonumber\\
\Theta_D(g;A)&=&\exp 2\pi i\Gamma_D(g;A).
\end{eqnarray}
 \(\Theta_D(g,A)\) satisfies the cocycle condition (2.10):
\[\Theta_D(g,\,A)\,\Theta_D(h,\,g\cdot A)=\Theta_D(gh,\,A),\qquad \mbox{for \(g,h\in D_0G\)}.\]
So if we define the action of \(D_0G\) on \({\cal A}^{\flat}(D)\times {\rm C}\)  by 
\[\left(g, (A,c)\right)\longrightarrow \left(g\cdot A\,,\Theta_D(g, A)\,c\right),\]
 we have an hermitian line bundle on \({\cal M}^{\flat}(S^3)\) with the transition function \(\Theta(g,A)\):
\begin{equation}
{\cal L}(S^3)={\cal A}^{\flat}(D)\times {\rm C}/D_0G \longrightarrow{\cal M}^{\flat}(S^3).
\end{equation}

Now we shall show that the line bundle \({\cal L}(S^3)\) has a connection with 
curvature given by the symplectic form \(-i\,\omega^{\flat}\).

We define the following 1-form \(\theta\) on \({\cal A}^{\flat}(D)\) :
\begin{equation}
\theta_A(a)= \frac{i}{48\pi^3}\int_D\,Tr[A^3\,a]+\int_D\,Tr\,\gamma^0_A(a)\,dx,\qquad 
\mbox{ for \(a\in T_A{\cal A}^{\flat}\)}.
\end{equation}
We note that \(\gamma^0\) is invariant under the action of \(D_0G\).
\(\theta\)  is a 0-cochain on the Lie group \(D_0G\) taking its value in the space of 1-forms on \({\cal A}^{\flat}(D)\), the coboundary  of \(\theta\) becomes
\begin{eqnarray*}
(\delta \,\theta_A(a))(g)&=& \frac{i}{48\pi^3}\int_D\,Tr[(g^{-1}Ag+g^{-1}dg)^3g^{-1}ag-A^3a]\\
&=& \widetilde d\,\Gamma(g,A)(a)+\frac{i}{48\pi^3}\int_D\,Tr[(AV^2+V^2A+2AVA)a],
\end{eqnarray*}
for \(g\in D_0G\),   
where \(V=dgg^{-1}\).   
But, since \(g\in D_0G\), we have \(V\vert S^3=0\), and 
\[\int_D\,Tr[(AV^2+V^2A+2AVA)a]=\int_{S^3}\,Tr[VA-AV)a]=0,\]
for \(a\in T_A{\cal A}^{\flat}(D)\).
Therefore
\begin{equation}
(\delta\,\theta_A(a))(g)= \widetilde d\,\Gamma_D(g,A)(a)= \frac{1}{2\pi i}(\widetilde d\,\log \Theta_D(g,A))(a).
\end{equation}
Thus  \(\theta_A\) gives a connection on the line bundle \({\cal L}(S^3)\).   The curvature form  becomes;
\begin{eqnarray*}
(\widetilde d\,\theta)_A(a,b)&=&\frac{i}{24\pi^3}\int_DTr[A^2(ab-ba)]+\phi_A(a,b)\nonumber\\
& =&\frac{i}{24\pi^3}\int_{S^3}Tr[A(ab-ba)]+\phi_A(a,b)\nonumber\\
&=&
-i\, \omega^{\prime}_A(a,b)+\phi_A(a,b).
\end{eqnarray*}

From Corollary 1.9 it follows that  \(\gamma^0_A\) 
and \(\phi_A\)  vanishes over 
\({\cal M}^{\flat}(S^3)\), hence 
\begin{eqnarray}
\theta_A(a)&=&\frac{i}{48\pi^3}\int_D\,Tr[A^3\,a],\\
\nonumber\\
(\widetilde d\,\theta)_A(a,b)
&=&-\frac{i}{24\pi^3}\int_{S^3}Tr[A(ab-ba)]=-i \,\omega^{\prime}_A(a,b),
\end{eqnarray}
for \(a,\,b\in T_A{\cal M}^{\flat}(S^3)\).

Thus we have proved the following
\begin{thm}~There exists a prequantization of the moduli space  \(\left({\cal M}^{\flat}(S^3),\,\omega^{\flat}\right) \), that is, 
 there exists a Hermitian line bundle  with connection \({\cal L}(S^3)\longrightarrow {\cal M}^{\flat}(S^3)\),   whose curvature is equal to the symplectic form  
\(-i\,\omega^{\flat}\).
\end{thm}

Similarly there is a hermitian line bundle with connection \({\cal L}((S^3)^{\prime})\) over the symplectic space \({\cal M}^{\flat}(S^3)\).   
\begin{equation}
{\cal L}((S^3)^{\prime})={\cal A}^{\flat}(D')\times {\rm C}/D_0^{\prime}G \longrightarrow{\cal M}^{\flat}(S^3).
\end{equation}
Remember that we identify  \({\cal M}^{\flat}(S^3)\) and  \({\cal M}^{\flat}((S^3)^{\prime})\) by Lemma 1.10.  \({\cal L}((S^3)^{\prime})\)  is given by the transition function 
\[\Theta_{D^{\prime}}(g^{\prime},A^{\prime})=\exp\,2\pi i\Gamma_{D^{\prime}}(g^{\prime},A^{\prime})\]
for \(g^{\prime}\in D^{\prime}_0G\) and \(A^{\prime}\in{\cal A}^{\flat}(D^{\prime})\).   
Here any 
 \(g^{\prime}\in D^{\prime}_0G \) is considered to be extended from \(D^{\prime}\) to \(S^4\) by putting it equal to the identity transformation on \(D\).   

The connection on \({\cal L}((S^3)^{\prime})\) is given by
\begin{equation}
\theta^{\prime}_A(a)=\frac{i}{48\pi^3}\int_{D^{\prime}}\,Tr[A^3\,a],\qquad 
\mbox{ for \(a\in T_A{\cal M}^{\flat}\)},
\end{equation}
Then \[
(\widetilde d\,\theta^{\prime})_A(a,b)=- (\widetilde d\,\theta)_A (a, b)=i\, \omega_A(a,b).\]

In fact the prequantum line bundles  \({\cal L}((S^3)^{\prime})\) and 
\({\cal L}(S^3)\) are in duality.   

For  \(A\in{\cal A}(D)\) and \({\cal A}(D^{\prime})\) such that 
\(A\vert S^3=A^{\prime}\vert (S^3)^{\prime}\), we write by \(A\vee A^{\prime}\) the connection on \(S^4\) obtained by patching up the two connections along the boundary.    We have 
\[(g\vee g^{\prime})\cdot(A\vee A^{\prime})=(g\cdot A)\vee (g^{\prime}\cdot A^{\prime}).\]

Now the Polyakov-Wiegmann formula  ( Lemma 2.2 ) says
\begin{equation}
C_5(F\,G)=C_5(F)+C_5(G)+\gamma(F,G),\,\qquad mod{\,\rm Z},
\end{equation}
for \(F,\,G\in {\cal G}(S^4)\).   Since 
\[\gamma(F,G)=\Gamma(F^{-1}dF,G)-C_5(G),\]
we have
\[\Gamma_D(g, A)+\Gamma_{D^{\prime}}(g^{\prime}, A^{\prime})=\Gamma(g\vee g^{\prime}, A\vee A^{\prime}),\qquad mod\,{\rm Z},\]
 for   
\(g\in D_0G\) and \(g^{\prime}\in D^{\prime }_0G\).   
Here we have used the fact \(\gamma(g\vee 1^{\prime},1\vee g^{\prime})=0\).

Therefore
\begin{equation}
\Theta_D(g, A)\,\Theta_{D^{\prime}}(g^{\prime}, A^{\prime})=\Theta_{S^4}(g\vee g^{\prime}, A\vee A^{\prime}).
\end{equation}

By virtue of this formula we have the homomorphism of line bundles:
\begin{equation}
{\cal L}(S^3)\times {\cal L}((S^3)^{\prime})\longrightarrow {\cal L}(\emptyset).
\end{equation}
It is given by
\begin{equation}
\left( [A^{\prime},c^{\prime}],\,[A,c]\right)\longrightarrow [A\vee A^{\prime},c\vee c^{\prime}], 
\end{equation}
where \(c\vee c^{\prime}\) in the right hand side means that \(c\in \pi^{-1}(A)\) and \(c^{\prime}\in \pi^{-1}(A^{\prime})\).   
Composed with the homomorphism \({\cal L}(S^4)\longrightarrow {\rm C}\), ( 2.15 ), 
we have the duality of 
\({\cal L}(S^3)\) and \({\cal L}((S^3)^{\prime})\):
\begin{equation}
{\cal L}(S^3)\times {\cal L}((S^3)^{\prime})\longrightarrow {\rm C}.
\end{equation}

\subsubsection{Chern-Simons functional  over \(D\) and \(D'\) }

We shall introduce the Chern-Simons functional on \(D\) as a section of the pullback line bundle 
\(p^{\ast}{\cal L}((S^3)^{\prime})\) that is horizontal with respect to the connection \(\theta\).   

We defined in ( 1.26 ) the map \(p:\,{\cal A}^{\flat}(D)\longrightarrow {\cal M}^{\flat}((S^3)')\) that induces the identification of \({\cal M}^{\flat}(S^3)\) and   \({\cal M}^{\flat}((S^3)')\) .   Let \(A\in{\cal A}^{\flat}(D)\).   We consider 
\[\left(A^{\prime},\,\exp 2\pi i{\rm CS}_{S^4}(A\vee A^{\prime})\right)\in {\cal A}^{\flat}(D^{\prime})\times{\rm C},\]
 for  \(A^{\prime}\in{\cal A}^{\flat}(D^{\prime})\) such that \([A^{\prime}]=p(A)\).   
Let \(A^{\prime \prime}\) be another flat connection on \(D^{\prime}\) with \(A^{\prime \prime}\vert S^3=A\vert S^3\).   There is a \(g^{\prime}\in{\cal G}_0(D^{\prime})\) such that \(A^{\prime \prime}=g^{\prime}\cdot A^{\prime}\).   
 So we have \((1\vee g^{\prime})\cdot (A\vee A^{\prime})=A\vee A^{\prime \prime}\).
Since 
\[{\rm CS}_{S^4}(A\vee A^{\prime \prime})={\rm CS}_{S^4}(A\vee A^{\prime})+\Gamma_{D^{\prime}}(A^{\prime},g^{\prime}),\qquad mod.\,{\rm Z},\]
it holds that
\[\exp 2\pi i\,{\rm CS}_{S^4}(A\vee A^{\prime \prime})=\Theta_{D^{\prime}}(A^{\prime},g^{\prime})\cdot \exp 2\pi i\,{\rm CS} _{S^4}(A\vee A^{\prime})\,,\]
so we have a well defined map
\begin{equation}
\exp 2\pi i\,{\rm CS}_D\,: \, A\longrightarrow [A^{\prime},\,\exp 2\pi i\,{\rm CS}_{S^4}(A\vee A^{\prime})]\in {\cal L}((S^3)^{\prime})_{p(A)},
\end{equation}
that is, we have a section \(\exp 2\pi i\,{\rm CS}_D\) of \(p^{\ast}{\cal L}((S^3)^{\prime})\) over \({\cal A}^{\flat}(D)\).   We call \({\rm CS}_D\)  the {\it Chern-Simons functional} over \(D\).

Similarly we have the Chern-Simons functional over \( D^{\prime}\) that is defined by 
\begin{equation}
\exp 2\pi i\,{\rm CS}_{D ^{\prime}}(A^{\prime})=[A,\,\exp 2\pi i\,{\rm CS}_{S^4}(A\vee A^{\prime})]\in{\cal L}(S^3)_{p^{\prime}(A^{\prime})},
\end{equation}
which is a section of the pull back line bundle \((p^{\prime})^{\ast}{\cal L}(S^3)\) by \(p^{\prime}:\,{\cal A}^{\flat}(D^{\prime})\longrightarrow{\cal M}^{\flat}(S^3)\).

Let \(\theta\) be the connection on the line bundle \({\cal L}(S^3)\).   
Let \(\nabla \) be the covariant derivative that corresponds to the connection form \(\theta\) defined by
\begin{equation}
\nabla=\widetilde d-\,i\pi\theta.
\end{equation}

\begin{prop}~
The functional \(\exp 2\pi i\,{\rm CS}_{D ^{\prime}}\) is a horizontal section of the line bundle 
\((p^{\prime})^{ \ast}{\cal L}(S^3)\):
\begin{equation}
\nabla \exp 2\pi i\,{\rm CS}_{D'} =0,
\end{equation}
and it gives the trivialization 
\[ (p^{\prime})^{\ast}{\cal L}(S^3)\simeq {\cal A}^{\flat}(D')\times{\rm C}. \]
\end{prop}

In fact we have from Lemma 2.2
\begin{eqnarray*}
(\,\widetilde d\,\exp\,2\pi i\,{\rm CS}_{D'}\,)_A(a)&=&\exp 2\pi i\,{\rm CS}_{D'}(A)\,(\,\widetilde d \,{\rm CS}_{D'}\,)_A(a)\\
&=& i\pi\theta_A(a)\,\exp 2\pi i\,{\rm CS}_{D'}(A).
\end{eqnarray*}

Similarly, for the pullback line bundle  \(p^{\ast}{\cal L}((S^3)')\simeq{\cal A}^{\flat}(D)\times{\rm C}\),  \(\exp 2\pi i\,{\rm CS}_{D}\) is a non-vanishing  horizontal section over \({\cal A}^{\flat}(D)\).

\subsection{Prequantization over a conformally flat 4-manifold}

In the preceeding sections we have defined, for \(M=S^4\), \(D\) and \(D'\),  
\begin{enumerate}
\item
the symplectic structure  \(\left({\cal M}^{\flat}(\partial M), \omega^{\flat}\right)\) on the \({\cal G}_0\)-orbit space of irreducible flat connections 
in the trivial bundle \(M\times G\).
\item
the prequantum line bundle, that is, 
the hermitian line bundles 
\({\cal L}(\partial M)\) over \({\cal M}^{\flat}(\partial M)\) with connection whose curvature is equal to the symplectic form of  \({\cal M}^{\flat}(\partial M)\).    
\item
the Chern-Simons functional \({\rm CS}_M\), such that \(\exp 2\pi i\,CS_M\) is a non-vanishing horizontal section 
 of the pull back line bundle of  \({\cal L}((\partial M)') \) by the boundary restriction map \(p:\,{\cal A}^{\flat}(M)\longrightarrow {\cal M}^{\flat}(\partial M)
\).   Here the upper prime indicates the opposite orientation.
\end{enumerate}
In the sequel we shall extend functorially these concepts to more general four-manifolds.

First we must make precise what the {\it Chern-Simons prequantization } means for general four-manifolds.   We follow the axiomatic approach proposed for the definition of Topological Field Theory (TFT) by Atiyah~\cite{A}, and likewise for the Wess-Zumino-Witten model treated in our previous work~\cite{K}.   

Let \({\cal X}_4\) be the conformal equivalence classes of all compact
conformally flat four-manifolds with boundary.     

Let \({\cal X}\) be the category whose objects are 
three-dimensional manifolds \(\Gamma\) which are disjoint union of round
\(S^3\)'s.   A morphism between three-dimensional manifolds
\(\Gamma_1\) and \(\Gamma_2\) is an oriented cobordism given by
\(\Sigma\in  {\cal X}_4\) with boundary   
\(\partial\Sigma=\Gamma_2\bigcup(\Gamma_1^{\prime})\), where
the upper prime indicates the opposite orientation.  

  Let 
\({\cal Y}\) be the category of complex line bundles.

To each \(M\in{\cal X}_4\) and \(\Gamma=\partial M\in{\cal X}\) we have the symplectic space 
\(( {\cal M}^{\flat}(\Gamma),\omega^{\flat} )\) and the 
boundary restriction map \(p:\,{\cal A}^{\flat}(M)\longrightarrow {\cal M}^{\flat}(\Gamma)\simeq {\cal M}^{\flat}(\Gamma')\).

 A Chern-Simons prequantization means a functor \(CS\) from the category
\({\cal X}\) to the category \({\cal Y}\) which assigns:
\begin{namelist}{axiom 1}
\item[{\bf CS1},]
 to each manifold \(\Gamma\in{\cal X}\), a complex line bundle
\(CS(\Gamma)\)  over the space \({\cal M}^{\flat}(\Gamma)\),
\item[{\bf CS2},]
to each \(M\in {\cal X}_4\), with  
\(\partial M=\Gamma\), a section 
\(\exp 2\pi i\, {\rm CS}_M\) over \({\cal A}^{\flat}(M)\) of the pullback line bundle \(p^{\ast} CS(\Gamma)  \) .
\end{namelist}

First the functor \(CS\) satisfies axioms of TFT.    The following axioms {\bf A1},{\bf A2} and
{\bf A3} represent in the category of line bundles the orientation reversal and 
the operation of disjoint union and contraction~\cite{A}.  

\begin{namelist}{ABC}
\item[{\bf A1}]( Involution ): 
 \begin{equation}
 CS(\Gamma^{\prime})= CS(\Gamma)^{\ast}
\end{equation}
 where \(\ast\) 
indicates the dual line bundle.
\item[{\bf A2}]( Multiplicativity ):   
\begin{equation}
 CS(\Gamma_1\cup\Gamma_2)= CS(\Gamma_1)\otimes
 CS(\Gamma_2).\end{equation}
\hskip 0.5cm
\item[{\bf A3}]( Associativity ):

For a composite cobordism 
 \(M=M_1\cup_{\Gamma_3}M_2\) such that 
 \(\partial M_1=\Gamma_1^{\prime}\cup\Gamma_3\) and 
 \(\partial M_2=\Gamma_2\cup\Gamma_3^{\prime}\), we have
\begin{equation}
\exp 2\pi i\,{\rm CS}_M(A)=<\exp 2\pi i\,{\rm CS}_{M_1}(A_1),
\exp 2\pi i\,{\rm CS}_{M_2}(A_2)>,\end{equation}
for any \(A\in{\cal A}^{\flat}(M)\), \(A_i=A\vert M_i\),
\(i=1,2\).     Where, if we let 
\(p_i:\,{\cal A}^{\flat}(M_i)\longrightarrow {\cal M}^{\flat}(\Gamma_i)\simeq  {\cal M}^{\flat}(\Gamma^{\prime}_i)\), \(i=1,2\), and \(p^i_3:\,{\cal A}^{\flat}(M_i)\longrightarrow {\cal M}^{\flat}(\Gamma_3)\simeq  {\cal M}^{\flat}(\Gamma^{\prime}_3)\),  \(i=1,2\), then 
\(<\,,\,>\) denotes the natural pairing
\begin{equation}
 {p_1}^{\ast}CS(\Gamma_1^{\prime})\otimes { p^1_3}^{\ast}CS(\Gamma_3)\otimes
  {p^2_3}^{\ast}CS(\Gamma_3^{\prime})
\otimes{p_2}^{\ast}CS(\Gamma_2)\longrightarrow   {p_1}^{\ast}CS(\Gamma_1^{\prime})\otimes
 {p_2}^{\ast}CS(\Gamma_2).
 \nonumber
 \end{equation}

\end{namelist}

From \({\bf A1}\) and \({\bf A2}\) it follows that any cobordism \(M\) between \(\Gamma_1\) and
\(\Gamma_2\) induces a homomorphism of sections of pullback line
bundles
\begin{equation}
\exp 2\pi i\,{\rm CS}_{M}:C^{\infty}(p_1^{\ast} CS(\Gamma_1))\longrightarrow
C^\infty ( p_2^{\ast} CS(\Gamma_2)).\end{equation}

We impose:
\begin{equation} CS(\phi)={\rm C}\quad\mbox{for \(\emptyset\) the empty
3-dimensional manifold.}
\end{equation}

The following axioms characterize the Chern-Simons prequantization.

\begin{namelist}{ABC}
\item[{\bf A4}] 
For each  
\(\Gamma\), \( CS(\Gamma)\) has a
connection with curavature   \(i\,\omega^{\flat}\).
\item[{\bf A5}]  
 \(\exp 2\pi i\,{\rm CS}_M\) is parallel with respect to the
induced connection on \(p^{\ast} CS(\Gamma)\) .
\end{namelist}

 Let  \(M\in{\cal X}_4\) and
 \(\Gamma=\bigcup_{i\in
I}\Gamma_i\) be the boundary of  \(M\).   Each oriented component \(\Gamma_i\) consists of  a
round \(S^3\), and is endowed with a parametrization
\(p_i:\,S^3\longrightarrow
\Gamma_i\).   We  
 distinguish positive  and negative parametrizations \(p_i
: S^3\longrightarrow \Gamma_i\,,
\,i\in I_{\pm}\), depending on whether \(p_i\) respects 
the orientation of  \(\Gamma_i\) or not.    
We assume further that 
the compact manifold  \(\widehat M\), that is obtained from \(M\) by capping the discs \(D_i\) to   \(\Gamma_i\) for \( i\in I_-\),  and \(D^{\prime}_i\) for \(i\in I_+\) respectively,  is connected and simply connected.    \(\widehat M\) itself is the boundary of a 5 dimensional manifold \(N\); \(\partial N=\widehat M\), since \(\widehat M\) is conformally flat.   

We define the line bundles \(CS(\emptyset)\), \(CS(S^3)\) and \(CS((S^3)')\) as follows:
\begin{eqnarray*}
CS(\emptyset)&=&{\cal L}(\emptyset), \\
CS(S^3)&=&{\cal L}((S^3)'), \\
 CS((S^3)')&=&{\cal L}(S^3).
 \end{eqnarray*}
The corresponding Chern-Simons functionals are 
\({\rm CS}_{S^4}\), \({\rm CS}_D\) and \({\rm CS}_{D'}\) respectively.    They satisfy the axioms {\bf A1}, {\bf A2}, {\bf A4} and {\bf A5}, for \(\Gamma=\emptyset, \, S^3\).

 We shall give in the following \(CS(\Gamma)\) and 
\({\rm CS}_M\).   

\(\widehat M\) being simply connected \({\cal M}^{ \flat}(\widehat M)\) is  one point, which we denote by \(\emptyset\).    The prequantum line bundle \({\cal L}(\emptyset)\) becomes the complex line \({\rm C}\).      Then the Chern-Simons functional on \(\widehat M\) is defined as in the case of \(S^4\).   
\begin{equation}
{\rm CS}_{\hat M}(A)=
\frac{i}{240\pi^3}\int_{N}\,Tr(A^5), \qquad \mbox{ mod.{\rm Z}}.
\end{equation}

For \(i\in I_-\oplus I_+\), the parametrization
defines the map
\(p_i:\,{\cal M}^{\flat}(\Gamma_i) \longrightarrow
{\cal M}^{\flat}(S^3)\).    Then we have the pull-back bundle
of
\({\cal L}(S^3)\) ( resp. \({\cal L}((S^3)^{\prime})\) ) by 
\(p_i\).   We define
\begin{eqnarray}
{\cal L}(\Gamma_i)&=& p_i^{\ast}{\cal L}(S^3)\quad\mbox{for
\(i\in I_-\),}\nonumber\\
{\cal L}(\Gamma_i)&=& p_i^{\ast}{\cal L}((S^3)^{\prime})\quad\mbox{for
\(i\in I_+\)},\end{eqnarray}
then we have respectively
\begin{eqnarray}
{\cal L}(\Gamma'_i)&=&p_i^{\ast}{\cal L}((S^3)^{\prime})\quad\mbox{for
\(i\in I_-\),}\nonumber\\
{\cal L}(\Gamma'_i)&=&p_i^{\ast}{\cal L}(S^3)
\quad\mbox{for
\(i\in I_+\)}.
\end{eqnarray}
Here the prime indicates the opposite orientation.

The line bundle \(CS(\Gamma)\)  
is defined by 
\begin{equation} 
CS(\Gamma)=\otimes_{i\in
I_-}{\cal L}(\Gamma_i)\otimes\otimes_{i\in I_+}
{\cal L}(\Gamma_i).
\end{equation}

Now let \(\alpha: S^3\longrightarrow S^3\) be the 
restriction on \(S^3\) of a 
conformal diffeomorphism on \(S^4\).   
First we suppose that \(\alpha\) preserves the orientation.  
Then, since the transition function \(\Theta_D\) is invariant
under \(\alpha\), the line bundle \({\cal L}(S^3)\) is invariant
under \(\alpha\).     If \(\alpha\) reverses the orientation
then \(D\) is mapped to \(D'\) and \(\Theta_D\) is changed to
\(\Theta_{D^{\prime}}\).   Then \(\alpha^{\ast}{\cal L}(S^3)\) becomes
\({\cal L}((S^3)^{\prime})\).    On the other hand the parametrizations 
\(p_i\) are uniquely defined up to composition with conformal 
diffeomorphisms.    Therefore \(CS(\Gamma)\) is well defined
for the conformal equivalence class of 
\(\Gamma\in{\cal X}\).

The dual of \(CS(\Gamma)\) is 
\begin{equation} 
CS(\Gamma^{\prime})= \otimes_{i\in
I_-}{\cal L}(\Gamma_i^{\prime})\otimes\otimes_{i\in I_+}
{\cal L}(\Gamma_i^{\prime}),\end{equation}
and the duality; \(CS(\Gamma)\times
CS(\Gamma')\longrightarrow{\rm C}\), is given from 
( 2.30 )
by:
\begin{eqnarray*}
&&<\,\otimes_{i\in I_-}[A'_i,c'_i]\otimes
\otimes_{i\in I_+}[B_i,d_i]\,,\,
\otimes_{i\in I_-}[A_i,c_i]\otimes
\otimes_{i\in I_+}[B'_i,d'_i]\,>\\ 
&&=\Pi_{i\in I_-}c_ic^{\prime}_i\cdot \Pi_{i\in I_+}
 d_id^{\prime}_i\cdot 
\exp\{-2\pi i\sum_{i\in I_-}{\rm CS}_{S^4}(A_i\vee A^{\prime}_i)
-2\pi i\sum_{i\in I_+}{\rm CS}_{S^4}(B_i\vee B^{\prime}_i)\}.
\end{eqnarray*}

\(CS(\Gamma)\) is seen to be a Hermitian line bundle with connection.   The connection \(\theta\) is given by the same formula as in (2.23):
\begin{equation}
\theta_A(a)=\,- \,\frac{i}{48\pi^3}\int_M\,Tr[A^3\,a],\qquad \mbox{ for \(a\in T_A{\cal M}^{\flat}(\Gamma)\)}.\end{equation}
The curvature of \(\theta\) is the symplectic form \(i\,\omega^{\flat}\) on \({\cal M}^{\flat}(\Gamma)\).

The Chern-Simons functional over \(M\) is defined as a section of the pull-back line bundle \(p^{\ast}CS(\Gamma)\) by the boundary restriction map \(p:{\cal A}^{\flat}(M)\longrightarrow {\cal M}^{\flat}(\Gamma)\).

We put 
\[\Psi=\otimes_{i\in I_-} \exp\,2\pi i\,{\rm CS}_{D_i}\otimes \otimes_{i\in I_+}\exp\,2\pi i\,{\rm CS}_{D_i^{\prime}}\,.\]
This is a section of the pullback line bundle of \(CS(\Gamma')\) by the boundary restriction map associated to \(\partial(\widehat M\setminus M)=\Gamma'\), here  \( \widehat M\setminus M=
\cup_{ i\in I_-}D_i\cup\cup_{i\in I_+}D_i^{\prime} \).     By the duality between \(CS(\Gamma)\) and  \(CS(\Gamma')\) there exists a section of the pull-back line bundle \(p^{\ast}CS(\Gamma)\).   This is by  definition the Chern-Simons functional \({\rm CS}_M\):
\begin{equation}
\left< \,\Psi(A')\,,\,\exp\,2\pi i\,{\rm CS}_M(A)\,\right>=\exp 2\pi i\,{\rm CS}_{\hat M}(A\vee A'),
\end{equation}
for \(A\in {\cal A}(M)\) and \(A'\in {\cal A}(\widehat M\setminus M)\).

From the construction the functor \(CS\) 
satisfies the axioms \((A1) - (A5)\).

\section{Wess-Zumino-Witten actions on prequantum line bundles}

\subsection{Four-dimensional WZW model}

We saw in 1.6 that the group of gauge transformations \({\cal G}(\partial M)\) acts on \({\cal M}^{\flat}(\partial M)\) by an infinitesimal symplectic automorphism.     
Now we investigate the problem of lifting this action to the line bundle \({\cal L}(\partial M)\).   For that we consider the Wess-Zumino-Witten line bundles and the associated abelian extesion of the group \({\cal G}(\partial M)\) that were introduced in~\cite{K,Mi}.     In 3.1 and 3.2 we shall review briefly the WZW model discussed in our previous work~\cite{K}.

Let \({\cal X}_4\),  \({\cal X}\) and  \({\cal Y}\) be as in the previous section.   
 WZW model is a functor \(WZ\) from \({\cal X}\) to \({\cal Y}\)  which
assigns:
\begin{namelist}{axiom 1}
\item[{\bf WZ1},]
 to each manifold $\Gamma\in{\cal X}$, a complex line bundle
$WZ(\Gamma)$  over the space $\Gamma G$,
\item[{\bf WZ2},]
to each $\Sigma\in {\cal X}_4$, with  
\(\partial\Sigma=\Gamma$, a non-vanishing section 
\( {\rm WZ}_{\Sigma}\) of the pullback line bundle \(r^{\ast} WZ(\Gamma)\) .
\end{namelist}
Where \(M G\) is the set of pointed smooth mappings from a manifold \(M\) to \(G\), and 
 \(r: \Sigma G\longrightarrow \Gamma G\) is the boundary restriction map.
 
 \(WZ\) is imposed to satisfy the TFT axioms,
\begin{namelist}{ABC}
\item[{\bf B1}]( Involution ): 
 \begin{equation}WZ(\Gamma^{\prime})=WZ(\Gamma)^{\ast}
\end{equation}
 where $\ast$ 
indicates the dual line bundle.
\item[{\bf B2}]( Multiplicativity ):   
\begin{equation}WZ(\Gamma_1\cup\Gamma_2)=WZ(\Gamma_1)\otimes
WZ(\Gamma_2).\end{equation}
\hskip 0.5cm
\item[{\bf B3}]( Associativity ):
For a composite cobordism 
 \(\Sigma=\Sigma_1\cup_{\Gamma_3}\Sigma_2\) such that 
 \(\partial\Sigma_1=\Gamma_1\cup\Gamma_3\) and 
 \(\partial\Sigma_2=\Gamma_2\cup\Gamma_3^{\prime}\), we have
\begin{equation}
{\rm WZ}_{\Sigma}(f)=<{\rm WZ}_{\Sigma_1}(f_1),
{\rm WZ}_{\Sigma_2}(f_2)>,\end{equation}
for any \(f\in\Sigma G\), \(f_i=f\vert\Sigma_i\),
i=1,2.
\end{namelist}
 \({\bf B1}\) and \({\bf B2}\) imply that 
any cobordism \(\Sigma \) between \(\Gamma_1\) and
\(\Gamma_2\) induces a homomorphism of sections of pullback line
bundles
\begin{equation}
WZ(\Sigma):C^{\infty}(\Sigma,r_1^{\ast}WZ(\Gamma_1))\longrightarrow
C^{\infty}(\Sigma,r_2^{\ast}WZ(\Gamma_2)).\end{equation}

 Besides these TFT axioms we demand that \(WZ\) satisfies; 
\begin{namelist}{ABC}
\item[{\bf B4}] For each 
\(\Gamma\),  \(WZ(\Gamma)\) has a
connection.
\item[{\bf B5}] 
\( {\rm WZ}_{\Sigma}\) is non-vanishing and horizontal with respect to the
induced connection on \(r^{\ast}WZ(\Gamma)\), 
\item[{\bf B6}] 
For each \(\Sigma\in {\cal X}_4\) with \(\Gamma=\partial\Sigma\),
on the pullback line bundle
\(r^{\ast}WZ(\Gamma)\) is defined a multiplication \(\ast\) with respect to which we
have  
\begin{equation}
{\rm WZ}_{\Sigma}(fg)={\rm WZ}_{\Sigma}(f)\ast {\rm WZ}_{\Sigma}(g)\qquad\mbox{for any $f,g\in\Sigma
G$}.
\end{equation}
\end{namelist}

We constructed in~\cite{K} the WZW model 
\[WZ:\, {\cal X}\,\ni ( \Gamma,\Sigma)\,\Longrightarrow \,(WZ(\Gamma),\,{\rm WZ}_{\Sigma}) \in {\cal Y}.\]
  There the construction of the functor \(WZ\) is performed by the same way as we mentioned in 2.6 for the Chern-Simons prequantization functor, that is, after we get \((\,WZ(S^3), \, {\rm WZ}_D\,)\), 
 \((\,WZ((S^3)'), \, {\rm WZ}_{D'}\,)\) and \((\,WZ( \emptyset ),\,{\rm WZ}_{\Sigma})\) for \({\Sigma}\in {\cal X}_4\) such that \(\partial {\Sigma}=\emptyset\), we investigate the relation between the orientation reversion and the duality for the line bundles \(WZ(S^3)\), etc., then we have  \((WZ(\Gamma),\,{\rm WZ}_{\Sigma})\) for \(\Gamma\in {\cal X}\) so that they satisfy the axioms.

 In the following we shall explain after~\cite{K}  how to construct  line bundles 
\(WZ(S^3)\), \(WZ(((S^3)')\) and \(WZ(\emptyset)\) for \(\emptyset=\partial S^4\), and the corresponding 
WZW functionals; \({\rm WZ}_D\), \({\rm WZ}_{D'}\) and \({\rm WZ}_{S^4}\).     We refer~\cite{K} for the proof.

\subsection{ WZW line bundles }

We define the action of \({\cal G}_0(D)=D_0G\) on  \(DG \times {\rm C}\) by 
\begin{equation}
h\cdot\,(f\,,c)=(\,fh\,,\,c\,\Theta_{D}(h\,,\,f^{ -1}df\,)\,,
 \end{equation}
for \(h\in D_0G\) and 
\((f,c)\in DG\times {\rm C}\).

Remember that, for  \(A=f^{ -1} df  \in{\cal A}(D)\) and a 
\(h\in D_0G\), 
\begin{eqnarray*}
\Theta_{D}(h,A)&=& \exp 2\pi i\,\Gamma_{D}(h\,,\,A).\\
\\
\Gamma_{D}(h\,,\,A)&=&
\frac{i}{24\pi^3}\int_{D} c^{1,1} (h,\,A)
+C_5(h\vee 1')\\
&=&\frac{i}{24\pi^3}\int_{D} c^{ 2,1} (f,\,h)
+C_5(h\vee 1').
\end{eqnarray*}

By definition we put
\begin{equation}
WZ((S^3)^{\prime})= DG\times {\rm C}/D_0G  
\end{equation}
It is a complex line bundle  over \(\Omega^3_0G\). 
We denote the equivalence class of \((f, c)\) by \([ f, c]\), and
 define the projection \[\pi:\,WZ((S^3)^{\prime}) \longrightarrow \Omega^3_0G\] by 
\(\pi([ f, c ])=  f\vert S^3\).    The transition function \(\chi(f,g)\) for 
\(f,g\in DG\) such that \(f\vert S^3=g\vert S^3\in \Omega^3_0G\) is given by  
\begin{equation}
\chi(f,g)=\Theta_{D}(f^{-1}g\,,\,f^{  -1} df).
\end{equation}

We have  another bundle on \(\Omega^3_0G\) .   
Let
\begin{equation}
WZ(S^3)= D'G\times {\rm C}/D^{\prime}_0G,
\end{equation}
where \(D^{\prime}_0G\) acts on \(D'G\times {\rm C}\) by 
\begin{equation}
h'\cdot\,(f'\,,\,c')=\left (f'h'\,,\,c'\Theta_{D'}(h',\,f^{ \prime  -1}  df')\right),
\end{equation}
for \(h'\in D^{\prime}_0G\) and 
\((f'\,,\,c'\,)\in D'G\times {\rm C}\).    
 The projection \(\pi:\,WZ(S^3) \longrightarrow \Omega^3_0G\) is given by 
 \([f',c' ]\longrightarrow f'|S^3\).   It is a complex line bundle.       The transition function \(\chi'(f',g')\) for 
\(f',g'\in D^{\prime}G\) such that \(f'\vert S^3=g'\vert S^3\) is given by 
 \begin{equation}
 \chi'(f',g')=\Theta_{D}(f^{\prime -1}g\,,\,f^{\prime  -1} df') .  
 \end{equation}
  Where
 \begin{eqnarray*}
\Theta_{D'}(h',A')&=& \exp 2\pi i\,\Gamma_{D'}(h',\,A').\\
\\
\Gamma_{D'}(h',\,A')&=&
\frac{i}{24\pi^3}\int_{D'} c^{1,1} (h',\,A')
+C_5(1\vee h')
\end{eqnarray*}

 \(WZ(S^3)\) and \(WZ((S^3)^{\prime}) \) are in duality .   To see that, let 
 \begin{equation}
WZ(\emptyset)=S^4G\times {\rm C}/S^4G
 \end{equation}
 with the action of \(S^4G\) on \(S^4G\times {\rm C}\) defined by 
 \(H\cdot (F,c)= (FH,\, c\,\Theta_{S^4}(H,\,F^{-1}dF)\,)\).   
 The duality is given by 
   \[WZ(S^3)\times WZ((S^3)^{\prime})\ni
 \left (\,[f',c']\,,[f , c]\,\right ) \longrightarrow \left[\,f\vee f'\,,\,cc'\,\exp 2\pi i\,C_5(f\vee f')\,\right]\in WZ(\emptyset),\]
 composed with the 
 isomorphism similar to (2.15):
 \[WZ(\emptyset) \,\ni [\,F, c\exp\{2\pi iC_5(F)\,]\}\longrightarrow \,c\in {\rm C}.\]

The WZW action functional \(WZ_{S^4}\) is defined by
\begin{equation}
WZ_{S^4}(f)=[f,\exp 2\pi i\,C_5(f)],\qquad\mbox{for \(f\in S^4G\)}.
\end{equation}

Let \(r:DG\longrightarrow S^3G\) and 
  \(r':D'G\longrightarrow (S^3)'G\) be the boundary restriction maps.   
The WZW action functional \({\rm WZ}_D\) is defined as a non-vanishing section of the pullback line bundle 
\(r^{\ast}WZ(S^3)\) .
\begin{equation}
{\rm WZ}_D(f)=[ f', \exp 2\pi i\, C_5(f\vee f')\,]\in WZ(S^3)\vert_{r(f)}, \qquad f\in DG.
\end{equation}
Similarly \({\rm WZ}_{D'}\) is defined as a section of \((r')^{\ast}WZ((S^3)')\).   
\begin{equation}
{\rm WZ}_{D'}(f')=[ f , \exp 2\pi i\, C_5(f\vee f')\,]\in WZ((S^3)')\vert_{r'(f')}, \qquad f'\in D'G.
\end{equation}

We have 
\begin{equation}
\langle {\rm WZ}_D(f),\,{\rm WZ}_{D'}(f')\rangle ={\rm  WZ}_{S^4}(f\vee f'), 
\end{equation}
for \(f\in DG\), \(f'\in D'G\) such that \(r(f)=r'(f')\).

\vspace{0.5cm}

The multiplication is defined in \(r^{\ast}WZ(S^3)\).   

For \(\lambda=[f',a' ]\in WZ(S^3)_{r(f)}\) and \(\mu=[g',b']\in WZ(S^3)_{r(g)}\), we put 
\begin{eqnarray}
(f,\lambda)\,\ast\,(g,\mu)&=&(f,\,[f',\,a'])\,\ast\, (g,\,[g',\,b'])\nonumber\\[0.2cm]
&=&(fg,\,[\,f'g'\,, a'b'
\exp\{2\pi i\gamma(f\vee f^{\prime},g\vee g^{\prime})]\,).
\end{eqnarray}

The right hand side does not depend on the choice of \([f',a']\) and  \([g',b']\) but on the equivalence classes \(\lambda\) and \(\mu\),
and the product is well defined.    

Then we have, from  Lemma 2.2, 
\begin{equation}
{\rm WZ}_D(f)\ast {\rm WZ}_D(g)={\rm WZ}_D(f g).
\end{equation}

Similarly we have the product  on \((r')^{\ast}WZ((S^3)')\) over \(D'G\).    We omit the product formula which is quite similar to the above.

\vspace{0.5cm}

Now we introduce the action of WZW action functionals on Chern-Simons functionals.   

Let \(\varphi \in r^{\ast}WZ(S^3)\)  and \(\alpha \in p^{\ast}CS(S^3)\).    Then  \(\varphi=(f, [f', a] ) \) with \(f\in DG\), \([f',a] \in WZ(S^3)_{r(f)}\), and \(\alpha=(A,[A',c])\) with \(A\in {\cal A}^{\flat}(D)\), \([A',c]\in CS(S^3)_{p(A)}\).   
The action of  \(r^{\ast}WZ(S^3)\) on \(p^{\ast}CS(S^3)\) is defined by
\begin{equation}
\varphi \ast \alpha =\left(\, f\cdot A, \, a c \,\exp 2\pi i\,\gamma(f\vee f', A\vee A')\,\right).
\end{equation}
Where, for \(F\in S^4G\) and \(B\in{\cal A}^{\flat}(S^4)\), we have put
\[\gamma(F,B)=\gamma(F,G)\qquad \,B=G^{-1}dG.\]

By the same argument as to prove ( 3.18 ) using Lemma 2.2 we can prove
the following
\begin{thm}   
\begin{enumerate}
\item
For \(f\in DG \) and \(A\in {\cal A}^{\flat}(D)\),
\begin{equation}
{\rm WZ}_D(f)\ast \exp 2\pi i\,{\rm CS}_D(A)=\exp 2\pi i\,{\rm CS}_D(f\cdot A).
\end{equation}
\item
For \(f'\in D'G \) and \(A' \in {\cal A}^{\flat}(D')\),
\begin{equation}
{\rm WZ}_{D'}(f')\ast \exp 2\pi i\,{\rm CS}_{D'}(A')=\exp 2\pi i\,{\rm CS}_{D'}(f' \cdot A').
\end{equation}
\end{enumerate}
\end{thm}

\vspace{0.5cm}

\subsection{ Abelian extensions of the group \(\Omega^3_0G\)}

The two-dimensional WZW action gives a central extension \(\widehat{LG}\) of the loop group \(L G\).   
   The total space of the
\(U(1)\)-principal bundle  \(\widehat{LG}\) was described as the set of
equivalence classes of pairs \((f,c)\in D^2G\times U(1)\),   where \(D^2\) is the
2-dimensional disc with boundary \(S^1\).    The equivalence relation was defined
on the basis of Polyakov-Wiegmann formula~\cite{PW}, as it was so in our
four-dimensional generalization treated above.

However the principal bundle associated to \(WZ(S^3)\) or \(WZ((S^3)^{\prime})\)  has not any natural group structure contrary 
to the case of the extension of loop group.      Instead 
J. Mickelsson gave an extension of \(\Omega^3_0G\) by 
the abelian group \(Map({\cal A}_3,U(1))\), where \({\cal A}_3\) is the space of connections on \(S^3\) .       
In the following we shall explain after~\cite{K,Mi, Mic} two extensions of 
  \(\Omega^3_0G\) by 
the abelian group \(Map({\cal A}_3,U(1))\) that are in duality.

We define the action of \(D_0G\) on  \(DG \times Map({\cal A}_3,U(1))\) by 
\begin{equation}
h\cdot\,(f\,,\lambda\,)=(\,fh\,,\,\lambda(\cdot) \Theta_{D}(h\,,\,f^{ -1}df\,)\,,
 \end{equation}
for \(h\in D_0G\) and 
\((f,\,\lambda\,)\in DG\times Map({\cal A}_3,U(1))\).

We consider the quotient space by this action;
\begin{equation}
\widehat{\Omega G}=DG\times Map({\cal A}_3,U(1))/D_0G.
\end{equation}
The equivalence class of \((f, \lambda)\) is denoted by \([ f, \lambda]\).   
The projection \(\pi: \widehat{\Omega G}\longrightarrow \Omega^3_0G\) is defined by \(\pi([f\,,\,\lambda])= f\vert S^3\).   Then \(\widehat{\Omega G}\) becomes a principal bundle over \(\Omega^3_0G\) with the structure group \(Map({\cal A}_3,U(1))\).   
The transition function is \(\chi(f,g)=\Theta_D(f^{-1}g, \,f^{-1}df )\) for \(f,\,g\in DG\) such that \(f\vert S^3=g\vert S^3\).   
Here the \(U(1)\) valued function \(\chi(\,f\,,g\,)\) is considered as a constant function in \(Map({\cal A}_3,U(1))\).

The group structure of  \(\widehat{\Omega G}\) is given by the Mickelsson's 
2-cocycle  on \(D\) which is defined by the following formula.
\begin{eqnarray}\gamma_D(f,g\,; A)&=&
\frac{i}{24\pi^3}\int_D(\delta c^{1,1})(\,f,\,g\,;\,A)
\nonumber\\[0.2cm]
&=& \frac{i}{24\pi^3}\int_{S^3}c^{2,0}(f,g\,;\,A)+\frac{i}{24\pi^3}\int_D\,c^{2,1}(f,g).
\end{eqnarray}  

We define the multiplication on \(DG\times Map({\cal A}_3,U(1))\) by
\begin{equation}
(f,\,\lambda)\bullet (g,\,\mu)=\left(fg\,
,\,\lambda(\cdot)\,\mu_{f}(\cdot)\exp\,2\pi i\,\gamma_{D}(\,f,g;\,A)\,\right),
\end{equation}
where \[
\mu_{f}(A)=\mu((f\vert S^3)^{-1}A(f\vert S^3)+(f\vert S^3)^{-1}d(f\vert S^3)).\]
Then \(DG\times Map({\cal A}_3,U(1))\) is endowed with a group structure.    
The associative law follows from (2.1) and (2.2);
\[\delta\,dc^{2,0}=d c^{3,0}=0,\qquad \delta\,c^{2,1}=0.\]
From the definition and Lemma 2.2 we can verify that  
\( [f,\lambda] = [g,\mu] \) if and only if there is a \(h\in D_0G\) such that \((g,\mu )=(f,\lambda)\bullet (h,\exp\,2\pi i\,C_5(h\vee1'))\), hence the set of elements  \((h,\,\exp 2\pi i\,C_5(h\vee 1')\,)\) with \(h\in D_0G\) forms a normal subgroup of \(DG\times Map({\cal A}_3,U(1))\).   Thus 
  \(\widehat{\Omega G}\) is endowed with the group structure.     
Since the group \(Map({\cal A}_3,U(1))\) is embedded as a normal subgroup of   \(\widehat{\Omega G}\),     \(\widehat{\Omega G}\)  is  an extension of \(\Omega^3_0G\) by the abelian group  \(Map({\cal A}_3,U(1))\).      

We have another extension of  \(\Omega^3_0G\) by  \(Map({\cal A}_3,U(1))\) .   
 We consider 
\begin{equation}
\widehat{\Omega' G}=D'G\times Map({\cal A}_3,U(1))/ D^{\prime}_0G,
\end{equation}
where the action of  \(D^{\prime}_0G\) is defined by
\begin{equation}
h' \cdot\,(f' \,,\lambda'\,)=(\,f' h' \,,\,\lambda' (\cdot) \Theta_{D'}(h'\,,\,f^{ \prime -1}df' \,)\,,
 \end{equation}
for \(h' \in D^{\prime}_0G\) and 
\((f' ,\,\lambda'\,)\in D'G\times Map({\cal A}_3,U(1))\).    The transition function is \(\chi'(f',g')\).   

The multiplication on \(\widehat{\Omega' G}\) is defined by the
 same way as above using the Mickelsson's 2-cocycle 
\(\gamma_{D'}(A;f',\,g')\) on \(D'\), and \(\widehat{\Omega' G}\) becomes a
extension of \(\Omega^3_0G\) by the abelian group 
 \(Map({\cal A}_3,U(1))\).   
 
 \vspace{0.5cm}
 
Let the group \(Map({\cal A}_3,U(1))\) act on \({\rm C}\) by
 \(\lambda\cdot c=\lambda(0)c\).   Then the associated line bundle to the principal bundle \(\widehat{\Omega' G}\) is  \(WZ(S^3)\), and the line bundle associated to \(\widehat{\Omega G}\) is \(WZ((S^3)')\).   

It is very easy to have a similar abelian extension of the group \({\cal G}(\Gamma)\) for \(\Gamma\in{\cal X}\).   In fact \(\Gamma\in {\cal X}\) being the disjoint union of \(S^3\), we have 
\[{\cal G}(\Gamma)=\otimes_{I_-\cup I_+}\Omega^3_0G.\]
Then the extension is given by the product bundle 
\begin{equation}
\widehat{{\cal G}(\Gamma)}=\otimes_{I_-}\widehat{\Omega G} \otimes \otimes_{I_+}
\widehat{\Omega' G}\,\longrightarrow \,{\cal G}(\Gamma),
\end{equation}
The multiplication (3.25) is extended to \(\widehat{{\cal G}(\Gamma)}\) by tensor products.

\subsection{WZW actions on prequantum line bundles}

We shall define  the action of  \(\widehat{\Omega G}\) on \({\cal L}(S^3)\).    

Let \( (\,f,\,\lambda\,)\in DG\times  Map({\cal A}_3,U(1)) \) and 
\((A,c)\in {\cal A}^{\flat}(D)\times{\rm C}\).   We put
\begin{equation}
\beta_D (f,A)=\frac{i}{24\pi^3}\int_D\,c^{1,1}(f,\,A).
\end{equation}
Note that the relation  \(\delta \beta_D=\gamma_D\) holds.   If \(f\in D_0G\) then  \(\Gamma_D(f,A)=\beta_D(f,A)+C_5(f\vee 1')\).    

The action of \( (\,f,\,\lambda\,) \) on \((A,c)\) is defined by
\begin{equation}
(\,f,\,\lambda\,)\bullet (A,c)=
\left(\,f\cdot A\,,\,c\lambda(A\vert S^3)\exp 2\pi i\,\beta_D(f\,,\,A)\,\right).
\end{equation}
It is a right action.   
The relation
\(\beta_D=\delta \gamma_D\) yields that it defines actually an action:
\[
(g,\mu)\bullet \left(\,(f,\lambda)\bullet (A,a)\right)=\left((f,\lambda)\bullet (g,\mu)\right)\,\bullet \,(A,a),\]
for \((f,\lambda),\,(g,\mu)\in DG\times Map({\cal A}_3,U(1))\) and \((A,a)\in {\cal A}^{\flat}(D)\times {\rm C}\).

From the definition of \({\cal L}(S^3)\) and the fact that \([h,\exp 2\pi i\,C_5(h\vee 1')]\) for \(h\in D_0G\) 
gives the unit of \(\widehat {\Omega G}\) we see that the above action descends to the action of 
 \(\widehat {\Omega' G}\) on \({\cal L}(S^3)\).   Thus we have proved the following theorem.  
 
\begin{thm}
The prequantum line bundle \({\cal L}(S^3)\) carries an action of  \(\widehat {\Omega G}\)  that is equivariant with respect to 
the infinitesimal symplectic action of \(\Omega^3_0G\) on the base space \({\cal M}^{\flat}(S^3)\).
\end{thm}

As for the reduction of \({\cal L}(S^3)\) by the action of \(\widehat{\Omega G}\) we have the following bundle over \({\cal N}^{\flat}(S^3)\).     
Let \[{\cal K}(S^3)=
{\cal A}^{\flat}(D)\times{\rm C}/\widehat{\Omega G},\]
and let \(\pi:\,{\cal K}(S^3)\longrightarrow {\cal N}^{\flat}(S^3)={\cal M}^{\flat}(S^3)/\Omega^3G \)
be the induced projection from \(\pi:\,{\cal L}(S^3)\longrightarrow {\cal M}^{\flat}(S^3)\).
Then \(\pi:\,{\cal K}(S^3)\longrightarrow {\cal N}^{\,\flat}(S^3)\) becomes a line bundle with the structure group \(Map({\cal A}_3, U(1))\).   
But \( {\cal N}^{\,\flat}(S^3)\) being one point we have 
\[ {\cal K}(S^3)\simeq {\rm C}.\]

Similarly there exists an action 
of  \(\widehat {\Omega' G}\) on  \({\cal L}((S^3)')\) that is equivariant with respect to 
the infinitesimal symplectic action \(\Omega^3_0G\) on \({\cal M}^{ \flat}((S^3)')\simeq {\cal M}^{\flat}(S^3)\).   In particular \(Map({\cal A}_3,U(1))\) acts on  \({\cal L}((S^3)')\), and 
 the reduction 
 \({\cal K}((S^3)')={\cal L}((S^3)')/Map({\cal A}_3,U(1))\) becomes a complex line \({\rm C}\).

We can easily generalize the above result to the prequantum line bundle \(CS(\Gamma)\) for \(\Gamma\in {\cal X}\) defined in (2.44).   The latter carries the equivariant action of \(\,\widehat{{\cal G}(\Gamma)}\) given by    the tensor product of the actions  
(3.30) defined on each factor.   

\begin{thm}
Let \(\Gamma\in {\cal X}\).   
The prequantum line bundle \(CS(\Gamma )\) carries an action of  \(\, \widehat{{\cal G}(\Gamma)}\)  that is equivariant with respect to 
the infinitesimal symplectic action of  \({{\cal G}(\Gamma)}\) on the base space \({\cal M}^{\flat}(\Gamma)\).   The reduction of this action becomes  
\[CS(\Gamma )/\widehat{{\cal G}(\Gamma)}\simeq {\rm C}^d,\]
where \(d\) is the number of the boundary components.
\end{thm}

\end{document}